\documentclass[11pt]{article}

\usepackage{color}
\usepackage{amsmath}
\usepackage{amssymb}
\usepackage{latexsym}

\newtheorem{assumption}{Assumption}

\def\qed{ \ \vrule width.2cm height.2cm depth0cm\smallskip}
\def \ep{\hbox{ }\hfill$\Box$}

\newenvironment{proof}{\noindent {\bf Proof.\/}}{$\qed$\vskip 0.1in}

\newcommand{\la}{\langle}
\newcommand{\ra}{\rangle}
\newcommand{\hP}{{\hat\dbP}}

\newcommand{\eps}{\varepsilon}

\newcommand{\ba}{\begin{array}}
\newcommand{\ea}{\end{array}}
\newcommand{\be}{\begin{equation}}
\newcommand{\ee}{\end{equation}}
\newcommand{\bea}{\begin{eqnarray}}
\newcommand{\eea}{\end{eqnarray}}
\newcommand{\beaa}{\begin{eqnarray*}}
\newcommand{\eeaa}{\end{eqnarray*}}

\def\dbE{\mathbb{E}}
\def\dbF{\mathbb{F}}

\def\dbH{\mathbb{H}}

\def\dbL{\mathbb{L}}

\def\dbN{\mathbb{N}}
\def\dbP{\mathbb{P}}
\def\dbR{\mathbb{R}}
\def\dbS{\mathbb{S}}
\def\dbQ{\mathbb{Q}}

\def\mrp{\mbox{\tiny M\!R\!P}}

\def\a{\alpha}
\def\b{\beta}

\def\e{\varepsilon}

\def\l{\lambda}

\def\si{\sigma}
\def\t{\tau}
\def\f{\varphi}
\def\th{\theta}
\def\o{\omega}

\def\O{\Omega}
\def\cA{{\cal A}}
\def\cB{{\cal B}}

\def\cD{{\cal D}}

\def\cF{{\cal F}}
\def\cG{{\cal G}}
\def\cH{{\cal H}}
\def\cI{{\cal I}}

\def\cM{{\cal M}}
\def\cN{{\cal N}}

\def\cP{{\cal P}}

\def\cT{{\cal T}}

\def\no{\noindent}

\def\ms{\medskip}
\def\bs{\bigskip}
\def\q{\quad}

\def\cd{\cdot}
\def\cds{\cdots}

\def\qed{ \hfill \vrule width.25cm height.25cm depth0cm\smallskip}

\newcommand{\basa}{\begin{assumption}}
\newcommand{\easa}{\end{assumption}}

\newcommand{\bas}{\begin{assum}}
\newcommand{\eas}{\end{assum}}

\def\limsup{\mathop{\overline{\rm lim}}}

\def\esup{\mathop{\rm ess\;sup}}

 \def\cd{\cdot}
\def\cds{\cdots}

\def\dis{\displaystyle}

\def\bx{\mathbf{x}}
\def\bQ{\mathbf{Q}}
\def\cad{{c\`{a}dl\`{a}g}}

\def\1{\mathbf{1}}

\def\:{\!:\!}
\def\reff#1{{\rm(\ref{#1})}}
\def \proof{{\noindent \it Proof.\quad}}

at 9pt

\begin{document}

\newtheorem{thm}{Theorem}[section]
\newtheorem{lem}[thm]{Lemma}
\newtheorem{cor}[thm]{Corollary}
\newtheorem{prop}[thm]{Proposition}
\newtheorem{rem}[thm]{Remark}
\newtheorem{eg}[thm]{Example}
\newtheorem{defn}[thm]{Definition}
\newtheorem{assum}[thm]{Assumption}

\renewcommand {\theequation}{\arabic{section}.\arabic{equation}}
\def\thesection{\arabic{section}}

\title{Quasi-sure Stochastic Analysis through Aggregation}

\date{Submitted: March 24, 2010. Accepted: August 19, 2011.}

\author{H.~Mete {\sc Soner}\footnote{ETH (Swiss Federal Institute of Technology),
Z\"urich and
Swiss Finance Institute, hmsoner@ethz.ch. Research partly supported by the
European Research Council under the grant 228053-FiRM.
Financial support from
the Swiss Finance Institute and the ETH Foundation
are also gratefully acknowledged.}
      \and Nizar {\sc Touzi}\footnote{CMAP, Ecole Polytechnique Paris,
      nizar.touzi@polytechnique.edu.
      Research supported by the Chair {\it Financial Risks} of the
      {\it Risk Foundation} sponsored by Soci\'et\'e
             G\'en\'erale, the Chair {\it Derivatives of the Future}
             sponsored by the {F\'ed\'eration Bancaire Fran\c{c}aise}, and
             the Chair {\it Finance and Sustainable Development}
             sponsored by EDF and Calyon. }
      \and Jianfeng {\sc Zhang}\footnote{University of
      Southern California, Department of Mathematics, jianfenz@usc.edu.
      Research supported in part by NSF grant DMS 06-31366 and DMS 10-08873.}
}\maketitle

\begin{abstract}
This paper is on developing stochastic analysis simultaneously under a
general family of probability measures
that are not dominated by a single probability measure.  The interest in this
question originates from
the probabilistic representations of fully nonlinear partial differential equations
and applications to mathematical finance. The existing literature relies either on the capacity
theory (Denis and Martini \cite{DM}), or on the underlying nonlinear partial differential equation
(Peng \cite{Peng-G}). In both approaches, the resulting theory requires certain smoothness, the so called quasi-sure continuity, of the
corresponding processes and random variables in terms of the underlying canonical process.
In this paper, we investigate this question for a larger class of ``non-smooth" processes,
but with a restricted family of non-dominated probability measures.
For smooth processes, our approach leads to similar results as in previous literature,
provided the restricted family satisfies an additional density property.

\vspace{5mm}

\noindent{\bf Key words:} non-dominated probability measures,
weak solutions of SDEs, uncertain volatility model, quasi-sure stochastic analysis.

\noindent{\bf AMS 2000 subject classifications:} 60H10, 60H30.
\end{abstract}
\newpage

\section{Introduction}
It is well known that all
probabilistic constructions crucially depend
on the underlying probability measure.
In particular,  all random variables and stochastic
processes are defined up to null sets of this
measure.  If, however, one needs to
develop stochastic analysis {\it simultaneously}
under a family of probability measures,
then  careful constructions are needed
as the null sets of different measures do not necessarily
coincide.  Of course,
when this family of measures is dominated
by a single measure this question trivializes
as  we can simply work with
the null sets of the dominating measure.
However, we are interested exactly in the cases where
there is no such dominating measure.
An interesting example of this situation
is provided in the study of financial markets
with uncertain volatility.  Then, essentially all measures
are orthogonal to each other.

Since for each probability measure we have a well
developed theory, for simultaneous
stochastic analysis, we are naturally led to the
following problem of aggregation.
Given a family of random variables or
stochastic processes, $X^\dbP$,  indexed by
probability measures $\dbP$, can one find an {\it aggregator}
$X$ that satisfies $X=X^\dbP$, $\dbP-$almost surely
for every probability measure $\dbP$?
This paper studies exactly
this abstract problem.
Once aggregation is achieved,
then essentially all classical results
of stochastic analysis generalize
as shown in Section \ref{sect-QS-Analysis} below.

This  probabilistic  question is also
closely related to  the theory
of second order backward stochastic differential equations
(2BSDE) introduced in \cite{cstv}.
These type of stochastic equations
have several applications in stochastic
optimal control, risk measures and
in the Markovian case, they provide
probabilistic representations for
fully nonlinear partial differential
equations.  A uniqueness result is also available
in the Markovian context
as proved in
\cite{cstv} using the theory of viscosity solutions.
Although the definition given in \cite{cstv}
does not require a special structure,
the non-Markovian case, however, is
better understood only recently.
Indeed, \cite{STZ09d} further develops the theory
and proves a general
existence and  uniqueness result by probabilistic techniques.
The aggregation result is a central
tool for this result and in
our accompanying papers  \cite{STZ09b,STZ09c,STZ09d}.
Our new approach to 2BSDE
is related to the quasi sure analysis introduced
by Denis and Martini \cite{DM} and the $G$-stochastic analysis of
Peng \cite{Peng-G}. These papers are motivated by the
volatility uncertainty in mathematical finance.
In such financial models the volatility of the underlying stock process
is only known to stay between two given bounds
$0\le \underline{a}<\overline{a}$.  Hence, in this context
one needs to define probabilistic objects simultaneously
for all probability measures under which the canonical process $B$
is a square integrable martingale with absolutely continuous quadratic
variation process satisfying
$$
\underline{a}dt\le d\la B\ra_t\le \overline{a}dt.
$$
Here $d\la B\ra_t$ is the quadratic variation process
of the canonical map $B$.  We denote the set of
all such measures by $\overline{\cP}_W$, but without
requiring the bounds $\underline{a}$ and $\overline{a}$,
see subsection \ref{ss.formulation}.

As argued above, stochastic analysis under a family
of measures naturally leads us to the
problem of aggregation.
This question, which is also outlined above,
is stated precisely in Section \ref{ss.aggregation},
Definition \ref{defn-aggregator}.
The main difficulty in aggregation originates from the fact that
the above family of probability measures are not dominated
by one single probability measure.  Hence the classical stochastic analysis
tools can not be applied simultaneously under all probability measures in this family.
As a specific example, let us consider the case of the stochastic integrals.
Given an appropriate integrand $H$,
the stochastic integral $I^\dbP_t=\int_0^t H_s dB_s$ can be defined
classically under each probability measure
$\dbP$.  However, these processes may depend on the
underlying probability measure.  On the other hand we
are free to redefine this integral outside the support
of $\dbP$.  So, if for example, we have two probability
measures $\dbP^1, \dbP^2$ that are orthogonal
to each other, see e.g. Example \ref{example-singular}, then the integrals are immediately aggregated
since the supports are disjoint.
However, for uncountably many probability measures,
conditions on $H$ or probability measures are needed.
Indeed,
in order to aggregate these integrals,
we need to construct a stochastic process
$I_t$ defined on all
of the probability space so that $I_t=I^\dbP_t$ for all $t$, $\dbP-$almost
surely.  Under smoothness assumptions on the integrand $H$ this
aggregation is possible
and a pointwise definition is provided by Karandikar \cite{Karandikar}
for c\`adl\`ag integrands $H$. Denis and Martini \cite{DM} uses the theory of
capacities and construct the integral for  {\it quasi-continuous} integrands,
as defined in that paper.
A different approach based on the underlying partial differential equation
was introduced by Peng \cite{Peng-G} yielding essentially the same results
as in \cite{DM}.  In Section
\ref{sect-QS-Analysis} below, we also provide a construction without any
restrictions on $H$ but in a slightly smaller class than $\overline{\cP}_W$.

For general stochastic processes or random variables,
an obvious consistency condition
(see Definition \ref{defn-consistency}, below)
is clearly needed for aggregation.
But  Example \ref{eg-no-aggregation} also shows that
this condition is  in general not sufficient.  So to obtain aggregation
under this minimal condition, we have two alternatives.
First is to
restrict the family of processes by requiring
smoothness.
Indeed the previous results of
Karandikar \cite{Karandikar}, Denis-Martini \cite{DM},
and Peng \cite{Peng-G} all belong to this case. A precise statement
is given in Section \ref{ss.aggregation} below.
The second approach is to slightly
restrict the class of non-dominated measures.
The main goal of this paper is to specify these restrictions
on the probability measures that allows us to prove aggregation
under only the consistency condition  \reff{consistency}.

Our main result, Theorem \ref{thm-aggregation},
is proved in Section \ref{sec-weak-aggregation}.
For this main aggregation result,
we assume that the class of probability measures
are constructed from a {\it separable} class of diffusion processes
as defined in subsection \ref{ss.separable}, Definition
\ref{defn-cA}.  This class of diffusion processes is somehow natural and
the conditions are motivated from stochastic optimal control.
Several simple examples
of such sets are also provided.
Indeed, the processes
obtained by a straightforward
concatenation of deterministic piece-wise constant
processes forms a separable class.
For most applications, this set would be sufficient.
However, we believe that working
with general separable class helps
our understanding of quasi-sure stochastic
analysis.

The construction of
a probability measure corresponding to a given diffusion process, however, contains
interesting technical details.
Indeed, given an $\dbF$-progressively measurable process $\a$,
we  would like to construct a
unique measure $\dbP^\a$.
For such a construction, we start with the Wiener measure $\dbP_0$
and assume that $\a$ takes values
in $\dbS^{>0}_d$
(symmetric, positive definite matrices)
and also satisfy $\int_0^t |\a_s| ds <\infty$ for all $t\ge 0$, $\dbP_0$-almost
surely.  We then consider the
$\dbP_0$ stochastic integral
\begin{equation}
\label{e.xa}
X^\a_t := \int_0^t \a^{{1\slash 2}}_s dB_s.
\end{equation}
Classically, the quadratic variation density of $X^\a$
under $\dbP_0$ is equal to $\a$.
We then set $\dbP^\a_S := \dbP_0 \circ (X^\a)^{-1}$
(here the subscript $S$ is for the strong formulation).
It is clear that $B$ under $\dbP^\a_S$ has the
same distribution as $X^\a$ under $\dbP_0$.
One can show that
the quadratic variation density of $B$ under
$\dbP^\a_S$  is equal to $a$
satisfying $a(X^\a(\omega))=\a(\omega)$
 (see Lemma \ref{l.technical} below for the existence of such $a$).
Hence, $\dbP^\a_S\in \overline{\cP}_W$.
Let $\overline{\cP}_S\subset \overline{\cP}_W$ be
the collection of all such local martingale measures $\dbP_S^\a$.
Barlow \cite{Barlow} has observed that  this
inclusion is strict.  Moreover,
this procedure changes the density of the quadratic variation
process to the above defined process $a$. Therefore to
be able to specify the quadratic variation a priori,
in subsection \ref{ss.local},
we consider the weak solutions of a
stochastic differential equation (\reff{SDE0}
below) which is closely related to
\reff{e.xa}.
This class of measures obtained as weak solutions
almost provides the necessary structure for aggregation.
The only additional structure we need is the
uniqueness of the map from the diffusion
process to the corresponding probability
measure.   Clearly, in general, there is
no uniqueness.  So we further restrict ourselves
into the class with uniqueness which we denote by
$\cA_W$.  This set and the probability measures generated by them,
$\cP_W$, are
defined in subsection \ref{ss.local}.

The implications of our aggregation result for quasi-sure stochastic analysis are
given in Section \ref{sect-QS-Analysis}.  In particular,
for a separable class of probability measures, we first construct
a quasi sure stochastic integral and then prove all classical results
such as Kolmogrov continuity criterion, martingale representation,
Ito's formula, Doob-Meyer decomposition and the Girsanov theorem.
All of them are proved as a straightforward application of
our main aggregation result.

If in addition the family of probability measures is dense in an appropriate sense,
then our aggregation approach provides the same result as the quasi-sure analysis.
These type of results, of course, require continuity
of all the maps in an appropriate sense.
The details of this approach
are investigated in our paper \cite{STZ09c}, see also Remark \ref{rem-dense}
in the context of the application to the hedging problem under uncertain volatility.
Notice that, in contrast with
\cite{DM}, our approach provides existence of an optimal hedging strategy,
but at the price of slightly restricting the family of probability measures.

The paper is organized as follows.  The local martingale measures
$\overline{\cP}_W$ and a universal filtration
are studied in Section \ref{sec-cPbar}.  The question
of aggregation is defined in Section \ref{ss.aggregation}.  In the next section,
we define $\cA_W$, $\cP_W$ and then the separable
class of diffusion processes.  The main aggregation result,
Theorem \ref{thm-aggregation}, is proved in
Section \ref{sec-weak-aggregation}.  The next section
generalizes several classical results of stochastic analysis
to the quasi-sure setting.  Section \ref{sec-UVM}
studies  the application to the hedging problem under uncertain volatility.
In Section \ref{ss.cPS} we  investigate the class $\overline\cP_S$
of mutually singular measures induced from strong formulation.
Finally, several examples concerning weak solutions and the proofs of several
technical results
are provided in the Appendix.
\vspace{15pt}

\noindent
{\bf Notations.}
We close this introduction with a list of notations introduced in the paper.
\begin{itemize}
\item {$\O:= \{\omega\in C(\dbR_+, \dbR^d):\o(0)=0\}$},
$B$ is the canonical process, $\dbP_0$ is the Wiener measure on $\Omega$.
\item For a given stochastic process $X$, $\dbF^X$ is the filtration generated by $X$.
\item $\dbF:= \dbF^B= \{\cF_t\}_{t\ge 0}$ is the filtration generated by $B$.
\item $\dbF^+:=\{\cF_t^+,t\ge 0\}$, where $\cF_t^+:=\cF_{t+}:=\bigcap_{s>t}\cF_s$,
\item $\cF^\dbP_t := \cF^+_t \vee \cN^\dbP(\cF^+_t)$
and $\overline\cF^\dbP_t := \cF^+_t \vee \cN^\dbP(\cF_\infty)$, where
$$\cN^\dbP(\cG) :=\left\{E\subset \O: ~\mbox{there exists}~ \tilde E\in\cG~~\mbox{such that}~~
E\subset\tilde E ~~\mbox{and}~~ \dbP[\tilde E]=0\right\}.$$
\item $\cN_\cP$ is the class of $\cP-$polar sets defined in Definition \ref{def-qs}.
\item $\hat\cF^\cP_t \;:=\;\bigcap_{\dbP\in\cP} \big(\cF^\dbP_t\vee \cN_\cP\big)$ is the universal
filtration defined in \reff{hatF}.
\item  $\cT$ is the set of all $\dbF-$stopping times
$\t$ taking values in $\dbR_+\cup\{\infty\}$.
\item  $\hat\cT^\cP$ is set of all $\hat\dbF^\cP-$stopping times.
\item $\la B\ra$ is the universally defined quadratic variation of $B$,
defined in subsection \ref{ss.formulation}.
\item $\hat a$ is the density of the quadratic variation $\la B \ra$, also
defined in subsection \ref{ss.formulation}.
\item $\dbS_d$ is the set of $d\times d$ symmetric matrices.
\item $\dbS_d^{>0}$ is the set of positive definite symmetric matrices.
\item $\overline{\cP}_W$ is the set of measures defined in subsection \ref{overlinecPW}.
\item $\overline{\cP}_S \subset \overline{\cP}_W$ is defined in the Introduction,
see also Lemma \ref{l.technical}.
\item $\overline{\cP}_{\mrp}\subset \overline{\cP}_W$
are the measures with the martingale representation
property, see \reff{overlinecPmrp}.
\item Sets $\cP_W$,  $\cP_S$, $\cP_{\mrp}$ are defined in subsection
\ref{ss.local} and section \ref{ss.cPS}, as the subsets of $\overline{\cP}_W$,  $\overline{\cP}_S$,
$\overline{\cP}_{\mrp}$ with the additional requirement of
weak uniqueness.
\item $\overline{\cA}$ is the set of integrable,
progressively measurable processes with values in $\dbS^{>0}_d$.
\item $\overline{\cA}_W := \bigcup_{\dbP\in\overline{\cP}_W} \overline{\cA}_W(\dbP)$ and
$\overline{\cA}_W(\dbP)$ is the set of diffusion matrices satisfying \reff{overlinecAW}.
\item ${\cA}_W$, $\cA_S$, $\cA_{\mrp}$ are defined as above using
$\cP_W$,  $\cP_S$, $\cP_{\mrp}$, see section \ref{ss.cPS}.
\item Sets $\O^{a}_{\hat\t}$, $\O^{a,b}_{\hat\t}$ and the stopping time $\th^{ab}$ are
defined in subsection \ref{ss.support}.
\item Function spaces $\dbL^{0}$, $\dbL^{p}(\dbP)$, $\hat\dbL^p$,
and the integrand spaces $\dbH^0$, $\dbH^{p}(\dbP^a)$,
$\dbH^{2}_{loc}(\dbP^a)$, $\hat\dbH^p$, $\hat\dbH^2_{loc}$
are defined in Section \ref{sect-QS-Analysis}.
\end{itemize}

\section{Non-dominated mutually singular probability measures}
\setcounter{equation}{0}
\label{sec-cPbar}

Let $\O:= C(\dbR_+, \dbR^d)$ be as above
and $\dbF=\dbF^B$ be the filtration generated by the canonical process
$B$.
Then it is well known that  this natural filtration
$\dbF$  is left-continuous, but is not right-continuous.
This paper makes use of the right-limiting filtration \
$\dbF^+$, the $\dbP-$completed filtration
$\dbF^\dbP := \{\cF^\dbP_t, t\ge 0\}$, and the
$\dbP-$augmented filtration
$\overline\dbF^\dbP:= \{\overline\cF^\dbP_t, t\ge 0\}$,
which are all right continuous.

\subsection{Local martingale measures}
\label{ss.formulation}

We say a probability measure $\dbP$ is a local martingale measure if
the canonical process $B$ is a local martingale under $\dbP$. It follows from
Karandikar \cite{Karandikar} that there exists an $\dbF-$progressively
measurable process, denoted as $\int_0^t B_s dB_s$, which coincides
with the It\^o's integral, $\dbP-$almost surely for all local martingale
measure $\dbP$. In particular, this provides a pathwise definition of
\beaa
\la B\ra_t := B_t B_t^T - 2 \int_0^t B_s dB_s
&\mbox{and}& \hat a_t:= \limsup_{\e\downarrow 0} {1\over \e}
[\la B\ra_t-\la B\ra_{t-\e}].
\eeaa
Clearly, $\la B\ra$ coincides with
the $\dbP-$quadratic variation of $B$,
$\dbP-$almost surely for all local martingale measure $\dbP$.

Let $\overline{\cP}_W$ denote the set
of all local martingale measures $\dbP$ such that
\bea
\label{overlinecPW}
\dbP\mbox{-almost surely},~
\la B\ra_t ~\mbox{is absolutely continuous in}~t ~\mbox{and}~
\hat a \ \mbox{takes values in}~ \dbS^{>0}_d,
\eea
where $\dbS^{>0}_d$ denotes the space of all
$d\times d$ real valued positive definite matrices.
We note that, for different $\dbP_1, \dbP_2\in\overline{\cP}_W$,
in general $\dbP_1$ and $\dbP_2$ are mutually singular,
as we see in the next simple example. Moreover, there is
no dominating measure for $\overline{\cP}_W$.

\begin{eg}
\label{example-singular}
{\rm Let $d=1$, $\dbP_1 := \dbP_0\circ (\sqrt{2}B)^{-1}$, and
$\O_i := \{\la B\ra_t = (1+i)t, t\ge 0\}$, $i=0, 1$. Then,
$\dbP_0, \dbP_1\in\overline{\cP}_W$,
$\dbP_0(\O_0) = \dbP_1(\O_1) = 1$, $\dbP_0(\O_1)=\dbP_1(\O_0)=0$, and $\O_0$ and $\O_1$ are disjoint.
That is, $\dbP_0$ and $\dbP_1$ are mutually singular.
\ep}
\end{eg}

In many applications, it is important that $\dbP\in\overline\cP_W$
has martingale representation property (MRP, for short),
i.e.~for any $(\overline\dbF^{\dbP},\dbP)$-local martingale $M$,
there exists a unique ($\dbP$-almost surely)
$\overline\dbF^{\dbP}$-progressively measurable
$\dbR^d$ valued process $H$  such that
\beaa
\int_0^t |\hat a^{1\slash 2}_s H_s|^2 ds<\infty ~~
\mbox{and}~~ M_t = M_0+\int_0^t H_s dB_s, ~~t\ge 0,
~~\dbP\mbox{-almost surely.}
\eeaa
We thus define
\bea
\label{overlinecPmrp}
\overline{\cP}_{\mrp} := \left\{\dbP\in\overline{\cP}_W:  B
~\mbox{has MRP under}~\dbP\right\}.
\eea
The inclusion $\overline{\cP}_{\mrp}\subset \overline{\cP}_W$ is strict
as shown in Example \ref{eg-MRP} below.

Another interesting subclass is the set $\overline\cP_S$
defined in the Introduction. Since in this paper it is
not directly used, we postpone its
discussion to Section \ref{ss.cPS}.

\subsection{A universal filtration}
\label{ss.filtration}

We now fix an arbitrary subset $\cP\subset \overline{\cP}_W$.
By a slight abuse of terminology, we define the following
notions introduced by Denis and Martini \cite{DM}.

\begin{defn}\label{def-qs}
{\rm (i)} We say that a property holds $\cP$-quasi-surely,
abbreviated as $\cP$-q.s., if it holds $\dbP$-almost surely for all $\dbP\in\cP$.
\\
{\rm (ii)} Denote $\cN_\cP := \cap_{\dbP\in\cP}\cN^\dbP(\cF_\infty)$
and we call $\cP$-polar sets the elements of $\cN_\cP$.
\\
{\rm (iii)} A probability measure $\dbP$ is called
absolutely continuous with respect to  $\cP$ if $\dbP(E)=0$
for all $E\in\cN_\cP$.
\end{defn}

In the stochastic analysis theory, it is usually assumed that the
filtered probability space satisfies the {\it usual hypotheses}.
However,  the key issue in the present paper is to develop stochastic
analysis tools simultaneously for non-dominated mutually singular
measures. In this case, we do not have a good filtration
satisfying the usual hypotheses under all the measures.
In this paper,  we shall use the following universal filtration
$\hat\dbF^{\cP}$ for
the mutually singular probability measures $\{\dbP, \dbP\in\cP\}$:
\bea\label{hatF}
\hat\dbF^{\cP} := \{\hat\cF^\cP_t\}_{t\ge 0}
&\mbox{where}& \hat\cF^\cP_t \;:=\;\bigcap_{\dbP\in\cP} \big(\cF^\dbP_t\vee \cN_\cP\big)
~\mbox{for}~
t\ge 0.
\eea

Moreover, we denote by $\cT$ (resp. $\hat\cT^\cP$) the
set of all $\dbF$-stopping times $\t$ (resp.,
 $\hat\dbF^\cP$-stopping times $\hat \t$) taking values in $\dbR_+\cup \{\infty\}$.

\begin{rem}
\label{rem-completion}
{\rm Notice that $\dbF^+ \subset\dbF^\dbP\subset \overline\dbF^{\dbP}$.
The reason for the choice of this completed filtration $\dbF^\dbP$ is as follows.
If we use the small filtration $\dbF^+$, then the crucial aggregation
result of Theorem \ref{thm-aggregation} below will not hold true.
On the other hand, if we use the augmented filtrations $\overline\dbF^\dbP$,
then Lemma \ref{lem-Pconsistent} below does not hold.
Consequently, in applications one will not be able to check the
consistency condition \reff{consistent} in Theorem \ref{thm-aggregation},
and thus will not be able to apply the aggregation result.
See also Remarks \ref{rem-augmentation0} and
 \ref{rem-augmentation} below.
However, this choice of the completed filtration does
 not cause any problems in the
applications.
\ep}
\end{rem}

We note that $\hat\dbF^\cP$ is right continuous and all $\cP$-polar
sets are contained in $\hat\cF^\cP_0$. But $\hat\dbF^\cP$ is not
complete under each $\dbP\in\cP$.  However, thanks to the
Lemma \ref{lem-version} below, all the properties
we need still hold under this filtration.

For any sub-$\si-$algebra $\cG$ of $\cF_\infty$ and any
probability measure $\dbP$, it is well-known that an
$\overline\cF_\infty^\dbP$-measurable random variable
$X$ is $[\cG\vee \cN^\dbP(\cF_\infty)]-$measurable if and
only if there exists a
$\cG$-measurable random variable $\tilde X$ such that $X=\tilde X$,
$\dbP$-almost surely.
The following result extends this property to processes and
states that one can always consider any process in its
$\dbF^+$-progressively
measurable version. Since $\dbF^+\subset\hat\dbF^\cP$, the
$\dbF^+$-progressively measurable version is also
$\hat\dbF^\cP$-progressively measurable. This important result
will be used throughout our analysis so as to consider any process in its
$\hat\dbF^\cP$-progressively measurable version.
However, we emphasize that the $\hat\dbF^\cP$-progressively
measurable version depends on the underlying probability measure $\dbP$.

\begin{lem}
\label{lem-version}
Let $\dbP$ be an arbitrary probability measure on the
canonical space $(\O,\cF_\infty)$, and let $X$ be an
$\overline \dbF^\dbP$-progressively measurable process.
Then, there exists a unique {\rm{(}}$\dbP$-almost surely{\rm{)}}
$\dbF^+$-progressively measurable process $\tilde X$ such that
$\tilde X=X$, $\dbP-$almost surely.
If, in addition, $X$ is c\`adl\`ag $\dbP$-almost surely,
then we can choose $\tilde X$ to be c\`adl\`ag $\dbP$-almost
surely.
\end{lem}
The proof is rather standard but it is provided in
Appendix for completeness. We note that, the identity $\tilde X = X$,
$\dbP$-almost surely, is equivalent to that they
are equal $dt\times d\dbP$-almost surely.
However, if both of them are \cad, then clearly
$\tilde X_t = X_t$, $0\le t\le 1$, $\dbP$-almost surely.

\section{Aggregation}
\label{ss.aggregation}

We are now in a position to define the problem.

\begin{defn}
\label{defn-aggregator}
{\rm{Let $\cP\subset\overline{\cP}_W$, and let $\{X^\dbP,\dbP\in\cP\}$
be a family of $\hat\dbF^\cP$-progressively measurable processes.
An $\hat\dbF^\cP$-progressively measurable process $X$
is called a}} $\cP$-{\it{aggregator}} {\rm{of the family
$\{X^\dbP,\dbP\in\cP\}$ if $X=X^\dbP$, $\dbP$-almost surely
for every $\dbP\in\cP$.}}
\end{defn}

Clearly, for any family $\{X^\dbP,\dbP\in\cP\}$ which can be aggregated,
the following consistency condition must hold.
\begin{defn}
\label{defn-consistency}
{\rm{We say that a family $\{X^\dbP,\dbP\in\cP\}$
satisfies the}} {\it{consistency condition}} {\rm{if, for any $\dbP_1, \dbP_2\in \cP$,
and $\hat\t\in \hat\cT^\cP$ satisfying $\dbP_1=\dbP_2$ on $\hat\cF^\cP_{\hat\t}$
we have}}
\bea
\label{consistency}
X^{\dbP_1} = X^{\dbP_2} ~\mbox{on}~[0, \hat\t], ~\dbP_1-\mbox{almost surely}.
\eea
\end{defn}

Example \ref{eg-no-aggregation} below shows that the
above condition is in general not sufficient. Therefore, we are left with
following two alternatives.
\begin{itemize}
\item{} Restrict the range of aggregating processes
by requiring that there exists a sequence of
$\hat\dbF^\cP$-progressively measurable processes
$\{X^n\}_{n\ge 1}$ such that $X^n\to X^\dbP$,
$\dbP$-almost surely as $n\to \infty$ for all $\dbP\in\cP$.
In this case, the $\cP$-aggregator is
$X := \limsup_{n\to\infty} X^n$. Moreover, the class
$\cP$ can be taken to be the largest possible class $\overline{\cP}_W$.
We observe that the aggregation results of Karandikar
\cite{Karandikar}, Denis-Martini \cite{DM}, and Peng \cite{Peng-G}
all belong to this case.
Under some regularity on the processes,
this condition holds.

\item{} Restrict the class $\cP$ of mutually
singular measures so that the consistency condition
\reff{consistency} is sufficient for the largest possible
family of processes $\{X^\dbP, \dbP\in\cP\}$.
This is the main goal of the present paper.
\end{itemize}

We close this section by constructing
an example in which the consistency
condition is not sufficient for aggregation.

\begin{eg}
\label{eg-no-aggregation}
{\rm Let $d=2$. First, for each
$x, y\in [1,2]$, let $\dbP^{x,y} := \dbP_0\circ (\sqrt{x} B^1, \sqrt{y} B^2)^{-1}$
and $\O_{x,y} := \{\la B^1\ra_t= xt, \la B^2\ra_t = yt, t\ge 0\}$.
Cleary for each $(x,y)$, $\dbP^{x,y}\in  \overline{\cP}_W$
and $\dbP^{x,y}[\O_{x,y}]=1$.
Next, for each $a\in [1,2]$, we define
\beaa
\dbP_a[E] := {1\over 2}\int_1^2 (\dbP^{a,z}[E] + \dbP^{z,a}[E])dz
&\mbox{for all}&
E\in \cF_\infty.
\eeaa
We claim that $\dbP_a\in \overline{\cP}_W$.
Indeed,  for any $t_1 < t_2$ and any  bounded
$\cF_{t_1}$-measurable random variable $\eta$,
we have
$$
2\dbE^{\dbP_a}[(B_{t_2}-B_{t_1})\eta]=
\int_1^2 \{\dbE^{\dbP^{a,z}}[(B_{t_2}-B_{t_1})\eta]+
\dbE^{\dbP^{z,a}}[(B_{t_2}-B_{t_1})\eta]\}dz=0.
$$
Hence $\dbP_a$ is a martingale measure. Similarly,
one can easily show that $I_2dt\le d\la B\ra_t \le 2I_2 dt$,
$\dbP_a$-almost surely,
where $I_2$ is the $2\times 2$ identity matrix.

For $a \in [1,2]$ set
$$
\O_a := \{\la B^1\ra_t = at, t\ge 0\}\cup\{\la B^2\ra_t = at, t\ge 0\}
\supseteq\cup_{z\in [1,2]}\left[\O_{a,z} \cup  \O_{z,a}\right]
$$
so that $\dbP_a[\O_a] = 1$.
Also for $a\neq b$, we have $\O_a\cap \O_b = \O_{a,b} \cup \O_{b,a}$
and thus $\dbP_a[\O_a\cap \O_b] = \dbP_b[\O_a\cap \O_b] = 0$.

Now let $\cP:=\{\dbP_a, a\in [1,2]\}$ and
set $X^a_t(\o) = a$ for all $t, \o$.
Notice that, for $a\neq b$, $\dbP_a$ and $\dbP_b$
disagree on $\cF^+_0 \subset \hat\cF^\cP_0$. Then
the consistency condition \reff{consistency} holds trivially.
However, we claim that there is no $\cP$-aggregator $X$
of the family $\{X^a, a\in [1,2]\}$. Indeed, if there is $X$ such that
$X = X^a$, $\dbP_a$-almost surely for all $a\in [1,2]$, then for any $a\in [1,2]$,
$$
1=\dbP_a[X^a_. = a] = \dbP_a[X_. = a]
= {1\over 2}\int_1^2\!\!\!\Big(\dbP^{a,z}[X_.=a] + \dbP^{z,a}[X_.=a]\Big)dz.
$$
Let $\l_n$ the Lebesgue measure on $[1,2]^n$ for integer $n\ge 1$.
Then, we have
\beaa
\l_1\Big(\{z: \dbP^{a,z}[X_.=a]=1\}\Big) =
\l_1\Big(\{z: \dbP^{z,a}[X_.=a]=1\}\Big) = 1,
~~\mbox{for all}~a\in [1,2].
\eeaa
Set $A_1:= \{(a,z): \dbP^{a,z}[X_.=a]=1\}$, $A_2:= \{(z,a): \dbP^{z,a}[X_.=a]=1\}$
so  that $\l_2(A_1) = \l_2(A_2)= 1$.  Moreover, $A_1\cap A_2 \subset \{(a,a): a\in (0,1]\}$
and  $\l_2(A_1\cap A_2) = 0$.  Now we directly calculate
that  $1\ge \l_2(A_1\cup A_2) = \l_2(A_1)+\l_2(A_2)-\l_2(A_1\cap A_2)=2$.
This contradiction implies that there is no aggregator.
\ep
}\rm
\end{eg}

\section{Separable classes of
mutually singular measures}
\setcounter{equation}{0}
\label{sec-separableclass}

The main goal of this section is to identify a condition on
the probability measures that yields aggregation
as defined in the previous section.
It is more convenient to specify this restriction through
the diffusion processes.  However, as we discussed in the
Introduction
there are technical difficulties
in the connection between the diffusion processes and
the probability measures.  Therefore, in the
first two subsections we will discuss the issue of
uniqueness of the mapping from
the diffusion process to
a martingale measure.
The separable class
of mutually singular measures are
defined in subsection \ref{ss.separable}
after a short discussion of the supports
of these measures in subsection \ref{ss.support}.

\subsection{Classes of diffusion matrices}
\label{ss.dm}

Let
\beaa
\overline{\cA}:= \Big\{a:\dbR_+ \to \dbS^{>0}_d
\ |\  \dbF \mbox{-progressively measurable and}
~\int_0^t |a_s|ds<\infty, ~\mbox{for all}~t\ge 0\Big\}.
\eeaa
For a given $\dbP\in \overline\cP_W$, let
\bea
\label{overlinecAW}
\overline{\cA}_W(\dbP) :=\Big\{a\in \overline\cA: a =
\hat a,~~\dbP\mbox{-almost surely}\Big\}.
\eea
Recall that $\hat a$ is the density of the quadratic variation of
$\la B\ra$ and is defined pointwise.
We also define
\beaa
\overline{\cA}_W := \bigcup_{\dbP\in\overline{\cP}_W} \overline{\cA}_W(\dbP).
\eeaa
A subtle technical point is that $\overline{\cA}_W$ is
strictly included in $\overline{\cA}$. In fact, the process
\beaa
a_t := \1_{\{\hat a_t \ge 2\}} + 3\1_{\{\hat a_t <2\}}
~~\mbox{is clearly in}~ \overline{\cA}~\backslash~ \overline{\cA}_W.
\eeaa

For any $\dbP\in \overline{\cP}_W$
and $a\in\overline{\cA}_W(\dbP)$, by the
L\'evy characterization,
the following It\^o's stochastic integral under $\dbP$ is a
$\dbP$-Brownian motion:
\bea
\label{WP}
W^{\dbP}_t := \int_0^t \hat a^{-{1\slash 2}}_sdB_s =
\int_0^t a^{-{1\slash 2}}_s dB_s, &t\ge 0.&
{ \dbP-\mbox{a.s.}}
\eea
Also since $B$ is the canonical process,
 $a = a(B_\cd)$ and thus
\be
\label{BWP}
dB_t= a^{{1\slash 2}}_t(B_\cd) dW^{\dbP}_t,~
\dbP\mbox{-almost surely,} ~\mbox{and}~W^{\dbP}_t
~\mbox{is a}~ \dbP\mbox{-Brownian motion}.
\ee

\subsection{Characterization by diffusion matrices}
\label{ss.local}

In view of \reff{BWP}, to construct a measure with a given
quadratic variation $a\in \overline{\cA}_W$,
we consider the stochastic differential equation,
\bea
\label{SDE0}
dX_t = a^{{1\slash 2}}_t(X_\cd) dB_t,~~\dbP_0\mbox{-almost surely}.
\eea
In this generality, we consider only weak
solutions $\dbP$ which we define next.
Although the following definition
is standard (see for example Stroock
\& Varadhan \cite{SV}),
we provide it for specificity.

\begin{defn}
\label{defn-weak} {\rm{Let $a$ be an element
of $\overline{\cA}_W$. \\
(i)~For $\dbF-$stopping times $\t_1 \le \t_2 \in\cT$
and a probability measure $\dbP^1$ on $\cF_{\t_1}$,
we say that $\dbP$ is a}} weak solution of \reff{SDE0}
on $[\t_1, \t_2]$ with initial condition $\dbP^1$, denoted as $\dbP\in \cP(\t_1, \t_2, \dbP^1, a)$, 
{\rm{if the followings hold:
 
1.~$\dbP=\dbP^1$ on $\cF_{\t_1}$ ;

2.~The canonical process $B_t$ is a
$\dbP$-local martingale on $[\t_1, \t_2]$;

3.~The process $W_t := \int_{\t_1}^t
a^{-{1\slash 2}}_s(B_\cd) dB_s$,
defined $\dbP-$almost surely for all
$t \in [\t_1, \t_2]$, is a
$\dbP$-Brownian Motion.
\\
(ii)~We say that the  equation \reff{SDE0}
has}} weak uniqueness on $[\t_1, \t_2]$
with initial condition $\dbP^1$
{\rm{if any two weak solutions
$\dbP$ and $\dbP'$ in $\cP(\t_1, \t_2, \dbP^1, a)$ satisfy
$\dbP=\dbP'$ on $\cF_{\t_2}$.
\\
(iii)~We say that  \reff{SDE0}}} has weak uniqueness
{\rm{if (ii) holds for any $\t_1, \t_2\in \cT$
and any initial condition
$\dbP^1$ on $\cF_{\t_1}$.}}
\end{defn}
We emphasize that the stopping times in this
definition are $\dbF$-stopping times.

Note that, for each $\dbP\in \overline\cP_W$ and  $a\in\overline\cA_W(\dbP)$,
$\dbP$ is a weak solution of \reff{SDE0} on $\dbR_+$ with initial value $\dbP(B_0=0) = 1$.
We also need uniqueness of this map
to characterize the
measure $\dbP$ in terms of the diffusion matrix $a$.
Indeed, if \reff{SDE0} with $a$ has weak uniqueness, we let
$\dbP^a\in \overline{\cP}_W$ be the unique
weak solution of \reff{SDE0} on $\dbR_+$
with initial condition $\dbP^a(B_0=0)=1$, and  define,
\bea
\label{cAW}
{\cA}_W :=\left\{ a\in \overline{\cA}_W:
~\mbox{\reff{SDE0} has weak uniqueness}\right\},
\qquad
\cP_W := \{\dbP^a: a\in \cA_W\}.
\eea
We also define
\bea
\label{cAMRP}
\cP_{\mrp} := \overline\cP_{\mrp}\cap \cP_W,\q
\cA_{\mrp} := \{a\in \cA_W: \dbP^a \in \cP_{\mrp}\}.
\eea
For notational simplicity, we denote
\bea
\label{Fa}
\dbF^a := \dbF^{\dbP^a},~~\overline\dbF^a :=
\overline\dbF^{\dbP^a}, &\mbox{for all}& a\in \cA_W.
\eea

It is clear that, for each $\dbP\in \cP_W$,
the weak uniqueness of the equation \reff{SDE0}
may depend on the version of
$a\in \overline{\cA}_W(\dbP)$.  This
is indeed the case and the following example
illustrates this observation.
\begin{eg}
\label{example-multipleweak}
{\rm Let $a_0(t) := 1$, $a_2(t):= 2$ and
\beaa
a_1(t) := 1 + \1_E \1_{(0,\infty)}(t),
&\mbox{where}& E:= \Big\{\limsup_{h\downarrow 0}{B_{h}
-B_0\over \sqrt{2h\ln\ln h^{-1}}} \neq 1\Big\}\in \cF^+_0.
\eeaa
Then clearly both $a_0$ and $a_2$ belong to $\cA_W$.
Also $a_1 = a_0$, $\dbP_0$-almost surely and
$a_1=a_2$, $\dbP^{a_2}$-almost surely.  Hence,
$a_1 \in \overline{\cA}_W(\dbP_0) \cap \overline{\cA}_W(\dbP^{a_2})$.
Therefore the equation \reff{SDE0} with coefficient $a_1$
has two weak solutions $\dbP_0$ and $\dbP^{a_2}$.
Thus $a_1\notin \cA_W$.
\ep

\begin{rem}
\label{rem-inverse}
{\rm  In this paper, we shall consider only those $\dbP\in \cP_W
\subset \overline{\cP}_W$.
However, we do not know whether this inclusion is
strict or not.  In other words, given an arbitrary
$\dbP\in \overline{\cP}_W$, can we always find one version
$a\in \overline{\cA}_W(\dbP)$ such that $a\in \cA_W$?
\ep}
\end{rem}
}\end{eg}

It is easy to construct
examples in $\cA_W$
in the Markovian context. Below, we
provide two classes of path dependent
diffusion processes in $\cA_W$.
These sets are in fact subsets of $\cA_S\subset \cA_W$,
which is defined in \reff{cPS} below.
We also construct some counter-examples  in
the Appendix. Denote
\bea
\label{Q}
\bQ:= \left\{ (t,\bx) \ :\ t \ge 0, \bx \in C([0,t], \dbR^d)
\right\}.
\eea

\begin{eg}
\label{eg-lip}
{\rm (Lipschitz coefficients)
Let
\beaa
a_t := \si^2(t, B_\cd) &\mbox{where}& \si: \bQ \to \dbS_d^{>0}
\eeaa
is Lebesgue measurable, uniformly Lipschitz continuous in
$\bx$ under the uniform norm, and
$\si^2(\cd,\mathbf{0})\in\overline{\cA}$.
Then  \reff{SDE0} has
a unique strong solution
and consequently $a\in\cA_W$.
\ep
}\end{eg}

\begin{eg}
\label{eg-Piecewise constant} {\rm
(Piecewise constant coefficients)
Let $a = \sum_{n=0}^\infty a_n \1_{[\t_n, \t_{n+1})}$
where $\{\t_n\}_{n\ge 0}\subset\cT$ is a nondecreasing
sequence of $\dbF-$stopping times with $\t_0=0$,
$\t_n\uparrow \infty$ as $n\to\infty$, and $a_n\in \cF_{\t_n}$
with values in $\dbS^{>0}_d$ for all $n$.
Again \reff{SDE0} has a unique strong solution and $a\in\cA_W$.

This example is in fact more involved than it looks like, mainly
due to the presence of the stopping times.
We relegate its proof to the Appendix.
\ep
}
\end{eg}

\subsection{Support of $\dbP^a$}
\label{ss.support}

In this subsection, we collect some properties of measures that are constructed in the previous subsection.
We fix a subset $\cA\subset \cA_W$, and
denote by $\cP:= \{\dbP^a: a\in\cA\}$ the corresponding
subset of $\cP_W$. In the sequel, we may also say
\beaa
\mbox{a property holds $\cA-$quasi surely if it holds $\cP-$quasi surely.}
\eeaa

For any $a\in \cA$ and any $\hat\dbF^\cP-$stopping time $\hat\t\in \hat\cT^\cP$, let
\bea
\label{Oa}
\O^{a}_{\hat\t}
:=
\bigcup_{n\ge 1}\Big\{\int_0^t \hat a_s ds =
\int_0^t a_sds,~\mbox{for all}~t\in [0,\hat\t+{1\over n}]\Big\}.
\eea
It is clear that
\bea
\label{Oabasic}
\O^a_{\hat\t}\in \hat\cF^\cP_{\hat\t},
~~\O^a_t~\mbox{is non-increasing in}~ t,
~~\O^a_{\hat\t+} = \O^a_{\hat\t}, ~~\mbox{and}
~ \dbP^a(\O^a_\infty) =1.
\eea
We next introduce the first disagreement
time of any $a, b\in\cA$, which plays a central role
in Section \ref{sec-weak-aggregation}:
\beaa
\th^{a,b} := \inf\Big\{t\ge 0: \int_0^t a_sds \neq \int_0^t b_sds\Big\},
\eeaa
and, for any $\hat\dbF^\cP-$stopping time $\hat\t\in\hat\cT^\cP$,
the agreement set of $a$ and $b$ up to $\hat\t$:
 \beaa
\O^{a,b}_{\hat\t}
:= \{\hat\t < \th^{a,b}\}\cup \{\hat\t=\th^{a,b}=\infty\}.
\eeaa
Here we use the convention that
$\inf\emptyset=\infty$.  It is obvious that
\bea
\label{Oabbasic}
\th^{a,b}\in \hat\cT^\cP,~~
\O^{a,b}_{\hat\t}\in \hat\cF^\cP_{\hat\t}
&\mbox{and}& \O^a_{\hat\t}\cap \O^b_{\hat\t}
\subset \O^{a,b}_{\hat\t}.
\eea

\begin{rem}
\label{rem-Abarsupport}
{\rm The above notations can be extended to all
diffusion processes $a, b\in \overline\cA$.
This will be important in Lemma \ref{lem-Astructure} below.
\ep}
\end{rem}

\subsection{Separability}
\label{ss.separable}

We are now in a position to
state the restrictions needed for the
main aggregation result Theorem \ref{thm-aggregation}.

\begin{defn}
\label{cA0}
{\rm A subset $\cA_0\subset {\cA}_W$ is called a}
generating class of diffusion coefficients {\rm if

{\rm (i)}~$\cA_0$ satisfies the concatenation property:
$
a \1_{[0,t)} + b \1_{[t,\infty)} \in \cA_0$ for $a, b\in\cA_0$, $t\ge 0$.

{\rm (ii)}~$\cA_0$ has constant disagreement times:
for all $a, b\in\cA_0$, $\th^{a,b}$ is a constant or,
equivalently,  $\O^{a,b}_t = \emptyset$ or $\O$ for all $t\ge 0$.}
\end{defn}

We note that the concatenation property is standard
in the stochastic control theory in order to establish
the dynamic programming principle, see, e.g.
page 5 in \cite{ST}. The constant disagreement times
property is important for both
Lemma \ref{lem-Pconsistent} below and
the aggregation result of Theorem \ref{thm-aggregation} below.
We will provide two examples of sets with these properties,
after stating the main restriction for the aggregation result.

\begin{defn}
\label{defn-cA}
{\rm We say $\cA$ is a}
separable class of diffusion coefficients generated by $\cA_0$
{\rm if $\cA_0\subset {\cA}_W$ is a
generating class of diffusion coefficients
and $\cA$ consists of all processes $a$  of the  form,
\bea
\label{Aa}
a = \sum_{n=0}^\infty \sum_{i=1}^\infty a^n_i
\1_{E^n_i}\1_{[\t_n, \t_{n+1})},
\eea
where $(a^n_i)_{i,n}\subset \cA_0$, $(\t_n)_n\subset\cT$ is
nondecreasing with $\t_0=0$  and

$\bullet$\quad $\inf\{n: \t_n=\infty\} <\infty$, $\t_n <\t_{n+1}$ whenever $\t_n<\infty$, 
and each $\t_n$ takes at most countably many values,

$\bullet$\quad for each $n$, $\{E^n_i, i\ge 1\}
\subset \cF_{\t_n}$ form a partition of $\O$.}
\end{defn}

We emphasize that  in
the previous definition the $\t_n$'s are $\dbF-$stopping times and $E^n_i \in \cF_{\t_n}$.
The following are two examples of
generating classes of diffusion coefficients.

\begin{eg}
\label{eg-deterministic}{\rm
Let $\cA_0\subset \overline{\cA}$ be the class of all deterministic
mappings. Then clearly $\cA_0\subset \cA_W$
and satisfies both properties (the concatenation
and the constant disagreement times properties)
of a generating class.
\ep}\end{eg}

\begin{eg}
\label{eg-Lipschitz}{\rm
Recall the set $\bQ$ defined in \reff{Q}. Let $\cD_0$
be a set of deterministic Lebesgue
measurable functions $\si: \bQ \to \dbS_d^{>0}$
satisfying,

-  $\si$ is uniformly Lipschitz continuous in $\bx$ under
$\dbL^\infty$-norm, and $\si^2(\cd,\mathbf{0})\in\overline{\cA}$ and

- for each $\bx\in C(\dbR_+,\dbR^d)$ and different $\si_1, \si_2\in \cD_0$,
the Lebesgue measure of the set $A(\si_1, \si_2, \bx)$ is equal to $0$,
where
\beaa
A(\si_1, \si_2,\bx) &:=& \Big\{t: \si_1(t,\bx|_{[0,t]}) = \si_2(t,\bx|_{[0,t]})\Big\}.
\eeaa
Let $\cD$ be the class of all possible concatenations
of $\cD_0$, i.e. $\si\in\cD$ takes the following form:
\beaa
\si(t,\bx) &:=& \sum_{i=0}^\infty \si_i(t,\bx)
\1_{[t_i, t_{i+1})}(t),~~(t,\bx)\in \bQ,
\eeaa
for some sequence $t_i\uparrow \infty$ and
$\si_i\in \cD_0$, $i\ge 0$. Let $\cA_0 := \{\si^2(t, B_\cd):
\si\in \cD\}$. It is immediate to check that
$\cA_0\subset{\cA}_W$ and
satisfies the concatenation  and the
constant disagreement times properties.
Thus it is also a generating class.
\ep}\end{eg}

We next prove several important
properties of separable classes.

\begin{prop}
\label{prop-cA}
Let $\cA$ be a separable class of
diffusion coefficients generated by $\cA_0$. Then
$\cA\subset {\cA}_W$, and
$\cA$-quasi surely is
equivalent to $\cA_0$-quasi surely.
Moreover, if $\cA_0\subset \cA_{\mrp}$, then $\cA\subset \cA_{\mrp}$.
\end{prop}

We need the following two lemmas
to prove this result. The first one provides
a convenient structure for the elements of $\cA$.

\begin{lem}
\label{lem-Astructure}
Let $\cA$ be a separable class of diffusion
coefficients generated by $\cA_0$. For any
$a\in\cA$ and $\dbF$-stopping time $\t \in \cT$,
there exist $\t \le \tilde\t \in\cT$, a sequence
$\{a_n, n\ge 1\}\subset \cA_0$, and a partition
$\{E_n, n\ge 1\}\subset \cF_{\t}$ of $\O$, such that
$\tilde\t>\t$ on $\{\t<\infty\}$
and
\beaa
a_t = \sum_{n\ge 1} a_n(t) \1_{E_n} &\mbox{for all}& t<\tilde\t.
\eeaa
In particular, $E_n\subset \O^{a,a_n}_\t$ and
consequently $\cup_n \O^{a,a_n}_\t =\O$. Moreover,
if $a$ takes the form \reff{Aa} and $\t\ge\t_n$,
then one can choose $\tilde\t \ge \t_{n+1}$.
\end{lem}
The proof of this lemma is straightforward, but
with technical notations.  Thus we postpone it to the Appendix.

We remark that at this point we do not know whether $a\in \cA_W$.
But the notations $\th^{a,a_n}$ and $\O^{a,a_n}_{\t}$ are well defined
as discussed in Remark \ref{rem-Abarsupport}. We recall from Definition \ref{defn-weak} that $\dbP\in \cP(\t_1, \t_2, \dbP^1, a)$ means $\dbP$ is a weak
solution of \reff{SDE0} on $[\tilde\t_1, \tilde\t_{2}]$
with coefficient $a$ and initial condition $\dbP^1$.

\begin{lem}
\label{lem-acombination}
Let $\t_1, \t_2\in\cT$ with $\t_1\le \t_2$,
$\{a^i, i\ge 1\}\subset \overline\cA_W$ {\rm (}not necessarily in ${\cA}_W${\rm )}
and $\{E_i, i\ge 1\}\subset\cF_{\t_1}$ be
a partition of $\O$. Let $\dbP^0$ be a probability measure on
$\cF_{\t_1}$ and  $\dbP^i\in \cP(\t_1, \t_2, \dbP^0, a^i)$  for $i\ge 1$.   Define
\beaa
\dbP(E) := \sum_{i\ge 1} \dbP^i(E\cap E_i) ~~\mbox{for all}~~ E\in\cF_{\t_2}\q\mbox{and}\q a_t := \sum_{i\ge 1} a^i_t \1_{E_i}, && t\in [\t_1, \t_2].
\eeaa
Then $\dbP\in \cP(\t_1, \t_2, \dbP^0, a)$.
\end{lem}

\proof Clearly, $\dbP=\dbP^0$ on $\cF_{\t_1}$.
It suffices to show that both $B_t$ and
$B_tB^T_t - \int_{\t_1}^t a_s ds$ are $\dbP$-local
martingales on $[\t_1, \t_2]$.

By a standard localization argument,
we may assume without loss of generality
that all the random variables  below are integrable.
Now for any $\t_1 \le \t_3 \le \t_4\le \t_2$ and
any bounded random variable $\eta \in \cF_{\t_3}$, we have
\beaa
\dbE^{\dbP}[(B_{\t_4}-B_{\t_3})\eta]& = &
\sum_{i\ge 1} \dbE^{\dbP^i}\Big[(B_{\t_4}-B_{\t_3})\eta \1_{E_i}\Big] \\
&= &\sum_{i\ge 1} \dbE^{\dbP^i}\Big[\dbE^{\dbP^i}
\Big(B_{\t_4}-B_{\t_3}|\cF_{\t_3}\Big)\eta \1_{E_i}\Big]
=0.
\eeaa
Therefore $B$ is a $\dbP$-local martingale on $[\t_1, \t_2]$.
Similarly one can show that
$B_tB^T_t - \int_{\t_1}^t a_s ds$ is also a $\dbP$-local
martingale on $[\t_1, \t_2]$.
\ep

\bs

\no{\it Proof of Proposition \ref{prop-cA}.}
Let $a\in \cA$ be given as in \reff{Aa}.

\ms

\no (i) We first show that $a\in\cA_W$. Fix $\th_1, \th_2\in\cT$ with
$\th_1\le \th_2$ and a probability measure
$\dbP^0$ on $\cF_{\th_1}$.  Set
\beaa
\tilde\t_0:=\th_1 &\mbox{and}& \tilde\t_n := (\t_n \vee\th_1)
\wedge\th_2, ~~n\ge 1.
\eeaa
We shall show that $\cP(\th_1, \th_2, \dbP^0, a)$ is a singleton, that is, the \reff{SDE0} on $[\th_1,\th_2]$ with
coefficient $a$ and initial condition $\dbP^0$ has a
unique weak solution.  To do this we prove 
by induction on $n$ that $\cP(\tilde\t_0,\tilde\t_n, \dbP^0, a)$ is a singleton.

First, let $n=1$.  We apply Lemma \ref{lem-Astructure}
with $\t=\tilde\t_0$ and choose $\tilde \t = \tilde\t_{1}$.
Then, $a_t = \sum_{i\ge 1} a_i(t) \1_{E_i}$
for all $t< \tilde\t_{1}$, where $a_i \in\cA_0$
and $\{E_i, i\ge 1\}\subset \cF_{\tilde\t_0}$ form a
partition of $\O$.
For $i\ge 1$, let $\dbP^{0,i}$ be the unique weak solution in $\cP(\tilde\t_0, \tilde\t_{1}, \dbP^0,a_i)$ and set
\beaa
\dbP^{0,a}(E) := \sum_{i\ge 1}
\dbP^{0,i}(E\cap E_i) &\mbox{for all}& E\in\cF_{\tilde\t_{1}}.
\eeaa
We use Lemma \ref{lem-acombination} to conclude that
$\dbP^{0,a}\in \cP(\tilde\t_0, \tilde\t_{1}, \dbP^0,a)$. On the other hand, suppose $\dbP\in \cP(\tilde\t_0, \tilde\t_{1}, \dbP^0,a)$ is an arbitrary weak solution. For each $i\ge 1$,
we define $\dbP^i$ by
\beaa
\dbP^i(E):= \dbP(E\cap E_i) + \dbP^{0,i}(E\cap (E_i)^c)
&\mbox{for all}& E\in\cF_{\tilde\t_1}.
\eeaa
We again use Lemma \ref{lem-acombination} and notice that $a \1_{E_i} + a_i \1_{(E_i)^c} = a_i$.
The result is that $\dbP^i\in  \cP(\tilde\t_0, \tilde\t_{1}, \dbP^0,a_i)$. Now by the uniqueness in $\cP(\tilde\t_0, \tilde\t_{1}, \dbP^0,a_i)$ we conclude that
$\dbP^i=\dbP^{0,i}$ on $\cF_{\tilde\t_1}$.
This , in turn,  implies that $\dbP(E\cap E_i)=
\dbP^{0,i}(E\cap E_i)$ for all $E\in\cF_{\tilde\t_1}$ and $i\ge 1$.
Therefore, $\dbP(E) = \sum_{i\ge 1} \dbP^{0,i}(E\cap E_i) =
 \dbP^{0,a}(E)$ for all $E\in \cF_{\tilde\t_1}$.
 Hence $\cP(\tilde\t_0,\tilde\t_1, \dbP^0, a)$ is a singleton.

We continue with the induction step.
Assume that $\cP(\tilde\t_0,\tilde\t_n, \dbP^0, a)$ is a singleton, and denote its unique element  by $\dbP^{n}$.
Without loss of generality, we assume
$\tilde\t_n <\tilde\t_{n+1}$.
Following the same arguments as above we know that
$\cP(\tilde\t_n,\tilde\t_{n+1}, \dbP^n, a)$ contains  a unique weak solution, denoted by $\dbP^{n+1}$.
Then both $B_t$ and $B_tB_t^T - \int_0^t a_s ds$
are $\dbP^{n+1}$-local martingales on $[\tilde\t_0, \tilde\t_n]$
and on $[\tilde\t_n, \tilde\t_{n+1}]$. This implies that $\dbP^{n+1}\in \cP(\tilde\t_0,\tilde\t_{n+1}, \dbP^0, a)$.
On the other hand, let $\dbP\in \cP(\tilde\t_0,\tilde\t_{n+1}, \dbP^0, a)$ be an arbitrary weak solution.
Since we also have  $\dbP\in  \cP(\tilde\t_0,\tilde\t_{n}, \dbP^0, a)$,
by the uniqueness in the induction assumption we
must have the equality $\dbP = \dbP^n$ on $\cF_{\tilde\t_n}$.
Therefore, $\dbP\in  \cP(\tilde\t_n,\tilde\t_{n+1}, \dbP^n, a)$.   Thus by uniqueness
$\dbP = \dbP^{n+1}$ on $\cF_{\tilde\t_{n+1}}$.
This proves the induction claim for $n+1$.

Finally, note that $\dbP^m(E)=\dbP^n(E)$ for all
$E\in\cF_{\tilde\t_n}$ and $m\ge n$.
Hence, we may define $\dbP^\infty(E):= \dbP^n(E)$
for $E\in \cF_{\tilde\t_n}$.  Since $\inf\{n: \t_n=\infty\} <\infty$,
then $\inf\{n: \tilde\t_n=\th_2\} <\infty$ and thus
$\cF_{\th_2} = \vee_{n\ge 1}\cF_{\tilde\t_n}$.
So we can uniquely extend $\dbP^\infty$ to $\cF_{\th_2}$.
Now we directly check that $\dbP^\infty\in \cP(\th_1,\th_2, \dbP^0, a)$ and is unique.

\ms
\no (ii) We next show that $\dbP^a(E) =0$ for all $\cA_0-$polar set $E$. Once again we
apply Lemma \ref{lem-Astructure} with $\t=\infty$.
Therefore $a_t = \sum_{i\ge 1} a_i(t)\1_{E_i}$
for all $t\ge 0$, where $\{a_i, i\ge 1\}\subset \cA_0$
and $\{E_i, i\ge 1\}\subset \cF_\infty$ form a
partition of $\O$. Now for any $\cA_0$-polar set $E$,
\beaa
\dbP^a(E) = \sum_{i\ge 1} \dbP^a(E\cap E_i) =
\sum_{i\ge 1} \dbP^{a_i}(E\cap E_i) = 0.
\eeaa
This clearly implies the equivalence between $\cA$-quasi surely and $\cA_0$-quasi surely.

\ms

\no (iii) We now assume $\cA_0\subset\cA_{\mrp}$ and show that $a\in \cA_{\mrp}$.
Let $M$ be a $\dbP^a$-local
martingale. We prove by induction on $n$ again that
$M$ has a martingale representation on $[0,\t_n]$ under $\dbP^a$
for each $n\ge 1$. This, together with the
assumption that $\inf\{n: \t_n=\infty\}<\infty$,
implies that $M$ has martingale representation
on $\dbR_+$ under $\dbP^a$, and thus proves that
$\dbP^a \in \cA_{\mrp}$.

Since $\t_0=0$, there is nothing to prove
in the case of $n=0$. Assume the result holds
on $[0, \t_n]$. Apply Lemma \ref{lem-Astructure}
with $\t=\t_n$ and recall that in this case
we can choose the $\tilde \t$ to be $\t_{n+1}$.
Hence $a_t = \sum_{i\ge 1} a_i(t) \1_{E_i}$, $t<\t_{n+1}$,
where $\{a_i, i\ge 1\}\subset \cA_0$ and
$\{E_i, i\ge 1\}\subset \cF_{\t_n}$ form a partition of $\O$.
For each $i\ge 1$, define
\beaa
M^i_t := [M_{t\wedge \t_{n+1}} - M_{\t_n}]\1_{E_i}\1_{[\t_n, \infty)}(t)
&\mbox{for all}& t\ge 0.
\eeaa
Then one can directly check that $M^i$ is a
$\dbP^{a_i}$-local martingale. Since $a_i\in\cA_0\subset \cA_{\mrp}$,
there exists $H^i$ such that $d M^i_t = H^i_t dB_t$,
$\dbP^{a_i}$-almost surely. Now define
$H_t := \sum_{i\ge 1} H^i_t \1_{E_i}$, $\t_n \le t<\t_{n+1}$.
Then we have $dM_t = H_t dB_t$,
$\t_n \le t <\t_{n+1}$, $\dbP^a$-almost surely.
\ep
\vspace{10pt}

We close this subsection by the following important example.
\begin{eg}
\label{eg-hP}
{\rm Assume $\cA_0$ consists of all
deterministic functions $a: \dbR_+\to\dbS^{>0}_d$
taking the form $a_t = \sum_{i=0}^{n-1}
a_{t_i} \1_{[t_i, t_{i+1})} + a_{t_n}\1_{[t_n,\infty)}$
where $t_i\in\dbQ$ and $a_{t_i}$ has rational entries.
This is a special case of Example \ref{eg-deterministic} and thus $\cA_0\subset\cA_W$.
In this case $\cA_0$ is countable.
Let $\cA_0 = \{a_i\}_{i\ge 1}$ and define
$\hP := \sum_{i=1}^\infty 2^{-i} \dbP^{a_i}$. Then
$\hP$ is a dominating probability measure of all $\dbP^a$, $a\in\cA$,
where $\cA$ is the separable class of diffusion
coefficients generated by $\cA_0$. Therefore, $\cA$-quasi surely
is equivalent to $\hP$-almost surely. Notice however that
$\cA$ is not countable.
\ep}
\end{eg}

\section{Quasi-sure aggregation}
\setcounter{equation}{0}
\label{sec-weak-aggregation}

In this section, we fix
\bea\label{Aseparable}
\mbox{ a separable class
$\cA$ of diffusion coefficients generated by $\cA_0$
}
\eea
and denote $\cP:=\{\dbP^a, a\in\cA\}$.
Then we prove the main aggregation result of this paper.

For this we recall that the notion of aggregation is defined
in Definition \ref{defn-aggregator} and
the notations $\th^{a,b}$ and $\O^{a,b}_{\hat\t}$ are introduced in
subsection \ref{ss.support}.

\begin{thm}[Quasi sure aggregation]\label{thm-aggregation}~ 
For $\cA$ satisfying \reff{Aseparable}, let $\{X^a, a\in\cA\}$
be a family of~~$\hat\dbF^\cP$-progressively measurable
processes. Then there exists a unique ($\cP-$q.s.)
$\cP$-aggregator $X$ if and only if $\{X^a, a\in\cA\}$ satisfies the consistency condition
\bea
\label{consistent}
X^a = X^{b},~\dbP^a-\mbox{almost surely on}~[0,\th^{a,b})
&\mbox{for any}~a\in\cA_0 ~\mbox{and}~ b\in\cA.
\eea

Moreover, if $X^a$ is c\`adl\`ag $\dbP^a$-almost
surely for all $a\in\cA$, then we can choose a $\cP$-q.s. {\cad}
version of the $\cP$-aggregator $X$.
\end{thm}

We note that the consistency condition \reff{consistent} is slightly
different from the condition \reff{consistency} before. The condition \reff{consistent}
is more natural in this framework and is more convenient to check in applications.
Before the proof of the theorem,
we first show that, for any $a, b\in\cA$,
the corresponding probability measures
$\dbP^a$ and $\dbP^b$ agree as long as $a$ and $b$ agree.

\begin{lem}
\label{lem-Pconsistent}
For $\cA$ satisfying \reff{Aseparable} and $a, b\in\cA$, $\th^{a,b}$ is an $\dbF$-stopping time taking countably many values and
\bea
\label{Pconsistent}
\dbP^a(E\cap \O^{a,b}_{\hat\t}) = \dbP^b(E\cap\O^{a,b}_{\hat\t})
&\mbox{for all}&
\hat\t\in\hat\cT^\cP~\mbox{and}~E\in\hat\cF^\cP_{\hat\t}.
\eea
\end{lem}

\proof (i) We first show that $\th^{a,b}$ is an
$\dbF$-stopping time. Fix an arbitrary time $t_0$. In view of Lemma \ref{lem-Astructure} with $\t=t_0$,
we assume without loss of generality that
\beaa
a_t = \sum_{n\ge 1} a_n(t) \1_{E_n} &\mbox{and}&
b_t = \sum_{n\ge 1} b_n(t) \1_{E_n} ~~
\mbox{for all}~~ t<\tilde \t,
\eeaa
where $\tilde\t>t_0$, $a_n, b_n \in\cA_0$
and  $\{E_n, n\ge 1\}\subset \cF_{t_0}$
form a partition of $\O$. Then
\beaa
\{\th^{a, b}\le t_0\} = \bigcup_n
\left[\{\th^{a_n, b_n}\le t_0\}\cap E_n\right].
\eeaa
By the constant disagreement times property
of $\cA_0$, $\th^{a_n,b_n}$ is a constant.
This implies that $\{\th^{a_n, b_n}\le t_0\}$
is equal to either $\emptyset$ or $\O$.
Since $E_n\in\cF_{t_0}$, we conclude that
$\{\th^{a, b}\le t_0\}\in\cF_{t_0}$ for all $t_0\ge 0$.
That is, $\th^{a,b}$ is an $\dbF$-stopping time.

\ms

\no (ii) We next show that $\th^{a,b}$
takes only countable many values.
In fact, by (i) we may now apply Lemma \ref{lem-Astructure}
with $\t=\th^{a,b}$. So we may write
\beaa
a_t = \sum_{n\ge 1} \tilde a_n(t)
\1_{\tilde E_n} &\mbox{and}& b_t =
\sum_{n\ge 1}\tilde b_n(t) \1_{\tilde E_n}
~~ \mbox{for all}~~ t<\tilde \th,
\eeaa
where $\tilde\th>\th^{a,b}$ or
$\tilde\th=\th^{a,b}=\infty$, $\tilde a_n, \tilde b_n \in\cA_0$, and
$\{\tilde E_n, n\ge 1\}\subset \cF_{\th^{a,b}}$ form a partition of $\O$.
Then it is clear that $\th^{a,b} = \th^{\tilde a_n, \tilde b_n}$
on $\tilde E_n$, for all $n\ge 1$. For each $n$,
by the constant disagreement times property of
$\cA_0$, $\th^{\tilde a_n,\tilde b_n}$ is  constant.
Hence $\th^{a,b}$ takes only countable many values.

\ms

\no (iii) We now prove \reff{Pconsistent}.
We first claim that,
\bea
\label{inclusion}
E\cap \O^{a,b}_{\hat\t} \in
\Big[\cF_{\th^{a,b}} \vee \cN^{\dbP^a}(\cF_\infty)\Big]
&\mbox{for any}& E\in \hat\cF^\cP_{\hat\t}.
\eea
Indeed, for any $t\ge 0$,
\beaa
E\cap \O^{a,b}_{\hat\t} \cap \{\th^{a,b}\le t\}
&=& E\cap \{\hat\t <\th^{a,b}\} \cap \{\th^{a,b}\le t\} \\
&=& \bigcup_{m\ge 1} \Big[E\cap \{\hat\t <\th^{a,b}\}
\cap \{\hat\t\le t-{1\over m}\}\cap\{\th^{a,b}\le t\}\Big].
\eeaa
By (i) above, $\{\th^{a,b}\le t\}\in \cF_t$.
For each $m\ge 1$,
\beaa
E\cap \{\hat\t <\th^{a,b}\}\cap \{\hat\t\le t-{1\over m}\}
\in \hat\cF^\cP_{t-{1\over m}} \subset \cF^+_{t-{1\over m}}
\vee \cN^{\dbP^a}(\cF_\infty) \subset \cF_t \vee \cN^{\dbP^a}(\cF_\infty),
\eeaa
and \reff{inclusion} follows.

By \reff{inclusion}, there exist $E^{a,i}, E^{b,i}\in\cF_{\th^{a,b}}$, $i=1,2$,
such that
\beaa
E^{a,1}\subset E\cap \O^{a,b}_{\hat\t}\subset E^{a,2},~ E^{b,1}\subset E\cap \O^{a,b}_{\hat\t}\subset E^{b,2}, ~\mbox{and}~ \dbP^a(E^{a,2}\backslash E^{a,1}) = \dbP^b(E^{b,2}\backslash E^{b,1}) = 0.
\eeaa
Define $E^1:= E^{a,1}\cup E^{b,1}$ and $E^2:= E^{a,2}
\cap E^{b,2}$, then
\beaa
E^1, E^2\in\cF_{\th^{a,b}},\q E^1 \subset E \subset E^2,~~\mbox{and}~~
\dbP^a(E^2\backslash E^1)=\dbP^b(E^2\backslash E^1)=0.
\eeaa
Thus $\dbP^a(E\cap \O^{a,b}_{\hat\t}) = \dbP^a(E^2)$ and
$\dbP^b(E\cap\O^{a,b}_{\hat\t}) = \dbP^b(E^2)$.
Finally, since $E^2\in \cF_{\th^{a,b}}$, following
the definition of $\dbP^a$ and $\dbP^b$, in particular
the uniqueness of weak solution of  \reff{SDE0} on
the interval $[0, \th^{a,b}]$, we conclude that
$\dbP^a(E^2)=\dbP^b(E^2)$. This implies \reff{Pconsistent} immediately.
\ep

\begin{rem}
\label{rem-augmentation0}
{\rm The property \reff{Pconsistent} is crucial for checking
 the consistency conditions in our aggregation result in
Theorem \ref{thm-aggregation}.  We note that \reff{Pconsistent}
does not hold if we replace the completed
$\si-$algebra $\cF^{a}_\t\cap \cF^{b}_\t$ with the augmented
$\si-$algebra $\overline\cF^{a}_\t\cap \overline\cF^{b}_\t$.
To see this, let $d=1$, $a_t := 1$, $b_t := 1+\1_{[1,\infty)}(t)$.
In this case, $\th^{a,b}=1$. Let $\t := 0$, $E := \O^a_1$.
One can easily check that $\O^{a,b}_0 = \O$, $\dbP^a(E) = 1$,
$\dbP^b(E)=0$. This implies that $E\in \overline\cF^{a}_0\cap \overline\cF^{b}_0$
and $E\subset \O^{a,b}_0$. However, $\dbP^a(E)=1\neq 0 = \dbP^b(E)$.
See also Remark \ref{rem-completion}.
\ep}
\end{rem}

\no
{\it Proof of Theorem \ref{thm-aggregation}.}
The uniqueness of $\cP-$aggregator is immediate.
By Lemma \ref{lem-Pconsistent} and the uniqueness of
weak solutions of \reff{SDE0} on $[0,\th^{a,b}]$,
we know $\dbP^a = \dbP^b$ on $\cF_{\th^{a,b}}$.
Then the existence of the $\cP$-aggregator obviously implies \reff{consistent}.
We now assume that the
condition \reff{consistent} holds and prove the existence of the $\cP$-aggregator.

We first  claim that, without loss of generality,
we may assume that $X^a$ is \cad. Indeed, suppose that
the theorem holds for {\cad} processes.  Then we construct a
$\cP$-aggregator for a family $\{X^a, a\in\cA\}$, not necessarily \cad,
as follows:

- If $|X^a|\le R$ for some constant $R>0$ and for all $a\in\cA$,
set $Y^a_t := \int_0^t X^a_s ds$. Then, the family
$\{Y^a, a\in\cA\}$ inherits the consistency
condition \reff{consistent}. Since $Y^a$ is continuous
for every $a\in\cA$, this family admits a
$\cP$-aggregator $Y$. Define $X_t := \limsup_{\e\to 0} {1\over \e}[Y_{t+\e}-Y_t]$.
Then one can verify directly that $X$ satisfies all the requirements.

- In the general case, set $X^{R,a}:= (-R)\vee X^a \wedge R$.
By the previous arguments there exists $\cP$-aggregator
$X^R$ of the family $\{X^{R,a}, a\in\cA\}$ and it is immediate
that $X := \limsup_{R\to \infty} X^R$ satisfies all the requirements.

We now assume that $X^a$ is \cad, $\dbP^a$-almost surely
for all $a\in\cA$. In this case, the consistency condition \reff{consistent}
is equivalent to
\bea
\label{consistent3}
X^a_t = X^{b}_t,~0\le t< \th^{a,b}, ~\dbP^a\mbox{-almost surely}
&\mbox{for any}~a\in\cA_0 ~\mbox{and}~ b\in\cA.
\eea
\no{\bf Step 1.}
We first introduce the following quotient sets of $\cA_0$.
For each $t$, and $a, b\in\cA_0$,
we say $a\sim_t b$ if $\O^{a,b}_{t}=\O$ (or, equivalently,
the constant disagreement time $\th^{a,b}\ge t$).
Then $\sim_t$ is an equivalence relationship in
$\cA_0$.  Thus one can form a partition of $\cA_0$
based on $\sim_t$. Pick an element from each partition set to
construct a quotient set $\cA_0(t)\subset\cA_0$.
That is, for any $a\in\cA_0$, there exists a unique
$b\in\cA_0(t)$ such that $\O^{a,b}_{t}=\O$. Recall the notation $\O^a_t$ defined in (\ref{Oa}).
By \reff{Oabbasic} and the constant disagreement
times property of $\cA_0$, we know that
$\{\O^a_{t}, a\in\cA_0(t)\}$ are disjoint.

\no{\bf Step 2.}  For fixed $t\in\dbR_+$, define
\bea
\label{xit}
\xi_t(\o) &:=& \sum_{a\in \cA_0(t)}
X^a_t(\o) \1_{\O^a_{t}}(\o)~~~\mbox{for all}~~~\o\in\O.
\eea
The above uncountable sum is well defined because
the sets $\{\O^a_{t}, a\in\cA_0(t)\}$ are disjoint.
In this step, we show that
\bea
\label{xi-aggregation}
\xi_t ~\mbox{is}~ \hat\cF^\cP_t\mbox{-measurable}
&\mbox{and}& \xi_t = X^a_t, ~~\dbP^a\mbox{-almost surely
for all}~ a\in\cA.
\eea
We prove this claim in the following three sub-cases.

\no{\bf 2.1.}
For each $a\in \cA_0(t)$, by definition $\xi_t = X^a_t$ on $\O^a_{t}$.
Equivalently $\{\xi_t\neq X^a_t\}\subset (\O^a_{t})^c$.
Moreover, by \reff{Oabasic}, $\dbP^a((\O^a_{t})^c)=0$.
Since $\O^a_{t}\in \cF^+_{t}$ and $\cF^a_{t}$
is complete under $\dbP^a$,
$\xi_t$ is $\cF^a_{t}$-measurable and $\dbP^a(\xi_t = X_t^a)=1$.

\no{\bf 2.2.}
Also, for each $a\in \cA_0$, there exists a
unique $b\in \cA_0(t)$ such that $a\sim_t b$.
Then $\xi_t = X_t^b$ on $\O^b_{t}$.
Since $\O_{t}^{a,b}=\O$, it follows
from Lemma \ref{lem-Pconsistent}
that $\dbP^a = \dbP^b$ on $\cF^+_{t}$ and
$\dbP^a(\O^b_{t}) = \dbP^b(\O^b_{t}) =1$.
Hence $\dbP^a(\xi_t = X_t^b)=1$.  Now by
the same argument as in the first case,
we can prove that $\xi_t$ is $\cF^a_{t}$-measurable.
Moreover, by the consistency condition \reff{consistent1},
$\dbP^a(X_t^a=X_t^b) = 1$.
This implies that $\dbP^a(\xi_t=X_t^a)=1$.

\no{\bf 2.3.}
Now consider  $a\in \cA$.  We apply
Lemma \ref{lem-Astructure} with $\t=t$.
This implies that  there exist a sequence
$\{a_j, j\ge 1\}\subset \cA_0$ such that
$\O = \cup_{j\ge 1} \O_{t}^{a,a_j}$. Then
\beaa
\{\xi_t \neq X_t^a\} &=&\bigcup_{j\ge 1}
\Big[\{\xi_t \neq X_t^a\} \cap \O_{t}^{a,a_j}\Big].
\eeaa
Now for each $j\ge 1$,
\beaa
\{\xi_t \neq X_t^a\} \cap \O_{t}^{a,a_j}
\subset \Big[\{\xi_t \neq X_t^{a_j}\}
\cap \O_{t}^{a,a_j}\Big] \bigcup
\Big[\{X_t^{a_j} \neq X_t^a\} \cap \O_{t}^{a,a_j}\Big].
\eeaa
Applying Lemma \ref{lem-Pconsistent}
and using the consistency condition \reff{consistent3},
we obtain
\beaa
\dbP^a\Big(\{X_t^{a_j} \neq X_t^a\}
\cap \O_{t}^{a,a_j}\Big)&=& \dbP^{a_j}
\Big(\{X_t^{a_j} \neq X_t^a\} \cap \O_{t}^{a,a_j}\Big)\\
&=& \dbP^{a_j}\Big(\{X_t^{a_j} \neq X_t^a\}
\cap \{t<\th^{a, a_j}\}\Big)=0.
\eeaa
Moreover, for $a_j\in\cA_0$, by the previous sub-case,
$\{\xi_t \neq X_t^{a_j}\} \in \cN^{\dbP^{a_j}}(\cF^+_{t})$.
Hence there exists $D\in\cF^+_{t}$ such that
$\dbP^{a_j}(D) = 0$ and $\{\xi_t \neq X_t^{a_j}\}\subset D$.
Therefore
\beaa
\{\xi_t \neq X_t^{a_j}\} \cap \O_{t}^{a,a_j}
\subset D \cap \O_{t}^{a,a_j} &\mbox{and}&
\dbP^a(D \cap \O_{t}^{a,a_j}) = \dbP^{a_j}(D \cap \O_{t}^{a,a_j}) = 0.
\eeaa
This means that $\{\xi_t \neq X_t^{a_j}\}
\cap \O_{t}^{a,a_j}\in \cN^{\dbP^a}(\cF^+_{t})$.
All of these together imply that
$\{\xi_t \neq X_t^a\}\in\cN^{\dbP^a}(\cF^+_{t})$.
Therefore, $\xi_t\in \cF^a_{t}$ and $\dbP^a(\xi_t=X_t^a)=1$.

Finally, since $\xi_t\in \cF^a_{t}$ for all
$a\in\cA$, we conclude that  $\xi_t\in\hat\cF^\cP_{t}$.
This completes the proof of \reff{xi-aggregation}.

\no{\bf Step 3.}
For each $n\ge 1$, set $t^n_i := {i\over n}, i\ge 0$ and
define
\beaa
X^{a,n} := X^a_0\1_{\{0\}}
+\sum_{i=1}^\infty X^a_{t^n_i} \1_{(t^n_{i-1}, t^n_{i}]}
~\mbox{for all}~a\in\cA &\mbox{and}
& X^n := \xi_0\1_{\{0\}}+\sum_{i=1}^\infty \xi_{t^n_i} \1_{(t^n_{i-1}, t^n_{i}]},
\eeaa
where $\xi_{t^n_i}$ is defined by \reff{xit}.
Let $\hat\dbF^n :=\{\hat\cF^\cP_{t+{1\over n}}, t\ge 0\}$.
By Step 2, $X^{a,n}, X^n$ are
$\hat\dbF^n$-progressively measurable
and $\dbP^a(X^n_t = X^{a,n}_t, t\ge 0) = 1$
for all $a\in\cA$. We now define
\beaa
X := \limsup_{n\to\infty} X^n.
\eeaa
Since $\hat\dbF^n$ is decreasing to $\hat\dbF^\cP$ and
$\hat\dbF^\cP$ is right continuous,
$X$ is $\hat\dbF^\cP$-progressively measurable.
Moreover, for each $a\in\cA$,
\beaa
\{X_t = X^a_t, t\ge 0\} \bigcap \{X ~\mbox{is \cad}\}
\supseteq \Big[\bigcap_{n\ge 1}\{X^{n}_t = X^{a,n}_t,
t\ge 0\}\Big] \bigcap \{X^a ~\mbox{is \cad}\}.
\eeaa
Therefore $X = X^a$ and $X$ is \cad,
$\dbP^a$-almost surely for all $a\in\cA$.
In particular, $X$ is \cad, $\cP$-quasi surely.
\ep

\ms

Let $\hat\t\in\hat\cT^\cP$ and $\{\xi^a, a\in \cA\}$ be
a family of $\hat\cF^\cP_{\hat\t}$-measurable
random variables. We say an  $\hat\cF^\cP_{\hat\t}$-measurable
random variable $\xi$ is a $\cP$-aggregator of the family $\{\xi^a, a\in \cA\}$ if
$\xi=\xi^a$, $\dbP^a$-almost surely for all $a\in\cA$. Note that
we may identify any $\hat\cF^\cP_{\hat\t}$-measurable
random variable $\xi$ with the $\hat\dbF^\cP$-progressively
measurable process $X_t := \xi \1_{[\hat\t, \infty)}$.
Then a direct consequence of Theorem  \ref{thm-aggregation} is
the following.

\begin{cor}
\label{cor-aggregation}
Let $\cA$ be satisfying \reff{Aseparable} and $\hat\t\in\hat\cT^\cP$. Then the family of $\hat\cF^\cP_{\hat\t}$-measurable
random variables $\{\xi^a, a\in \cA\}$ has a unique
{\rm(}$\cP$-q.s.{\rm)}
$\cP$-aggregator $\xi$ if and only if the
following consistency condition holds:
\be
\label{consistent1}
\xi^a = \xi^{b} ~\mbox{on}~\O^{a,b}_{\hat\t},~~
\dbP^a\mbox{-almost surely for any}~
a\in\cA_0 ~\mbox{and}~ b\in\cA.
\ee
\end{cor}
\vspace{5mm}

For the next result, we
recall that the $\dbP$-Brownian motion $W^\dbP$ is defined in \reff{WP}.
As a direct consequence of Theorem \ref{thm-aggregation},
the following result defines the $\cP$-Brownian motion.

\begin{cor}\label{cor-quadratic}
For $\cA$ satisfying \reff{Aseparable}, the family $\{W^{\dbP^a},a\in\cA\}$ admits a unique $\cP$-aggregator $W$.
Since $W^{\dbP^a}$ is a $\dbP^a$-Brownian motion
for every $a\in\cA$, we call $W$ a $\cP$-universal Brownian motion.
\end{cor}

\proof  Let $a, b \in \cA$. For each $n$, denote 
\beaa
\t_n := \inf\Big\{ t\ge 0: \int_0^t |\hat a_s|ds \ge n\Big\}\wedge \th^{a,b}.
\eeaa
Then $B_{\cd\wedge \t_n}$ is a $\dbP^b$-square integrable martingale. By standard construction of stochastic integral, see e.g. \cite{KS} Proposition 2.6, there exist $\dbF$-adapted simple processes $\b^{b,m}$ such that 
\bea
\label{betam}
\lim_{m\to\infty}\dbE^{\dbP^b}\Big\{\int_0^{\t_n} |\hat a_s^{1\over 2}(\b^{b,m}_s - \hat a_s^{-{1\over 2}})|^2 ds \Big\} =0.
\eea
Define the universal process
\beaa
W^{b,m}_t := \int_0^t \b^{b,m}_s dB_s.
\eeaa
Then
\bea
\label{Wm}
\lim_{m\to\infty}\dbE^{\dbP^b}\Big\{\sup_{0\le t\le \t_n}\Big|W^{b,m}_t- W^{\dbP^b}_t\Big|^2  \Big\} =0.
\eea
By Lemma \ref{lem-version}, all the processes in (\ref{betam}) and (\ref{Wm}) can be viewed as $\dbF$-adapted. Since $\t_n \le \th^{a,b}$, applying Lemma  \ref{lem-Pconsistent} we obtain from \reff{betam} and \reff{Wm} that
\beaa
 \lim_{m\to\infty}\dbE^{\dbP^a}\Big\{\int_0^{\t_n} |\hat a_s^{1\over 2}(\b^{b,m}_s - \hat a_s^{-{1\over 2}})|^2 ds \Big\} =0,\q \lim_{m\to\infty}\dbE^{\dbP^a}\Big\{\sup_{0\le t\le \t_n}\Big|W^{b,m}_t- W^{\dbP^b}_t\Big|^2  \Big\} =0.
\eeaa
The first limit above implies that 
\beaa
\lim_{m\to\infty}\dbE^{\dbP^a}\Big\{\sup_{0\le t\le \t_n}\Big|W^{b,m}_t- W^{\dbP^a}_t\Big|^2  \Big\} =0,
\eeaa
which, together with the second limit above, in turn leads to
\beaa
W^{\dbP^a}_t = W^{\dbP^b}_t, ~~ 0\le t\le \t_n,\q \dbP^a-\mbox{a.s.}
\eeaa
Clearly $\t_n \uparrow \th^{a,b}$ as $n\to\infty$. Then
\beaa
W^{\dbP^a}_t = W^{\dbP^b}_t, ~~ 0\le t<\th^{a, b},\q \dbP^a-\mbox{a.s.}
\eeaa
That is,  the family $\{W^{\dbP^a}, a\in\cA\}$
satisfies the consistency condition \reff{consistent}.
We then apply   Theorem \ref{thm-aggregation} directly
to obtain the $\cP-$aggregator $W$.
\ep

\bs

The $\cP-$Brownian motion $W$ is our first example of a
stochastic integral defined simultaneously under all
$\dbP^a$, $a\in\cA$:
\bea\label{W}
W_t &=& \int_0^t \hat a_s^{-1/2} dB_s, ~~t\ge 0,~~\cP-\mbox{q.s.}
\eea
We will investigate in detail the
universal integration
in Section \ref{sect-QS-Analysis}.

\begin{rem}
\label{rem-augmentation}
{\rm Although $a$ and $W^{\dbP^a}$ are
$\dbF$-progressively measurable,
from Theorem \ref{thm-aggregation}
we can only deduce that $\hat a$ and $W$ are
$\hat\dbF^\cP$-progressively measurable.
On the other hand, if we take a version of
$W^{\dbP^a}$ that is progressively measurable
to the augmented filtration $\overline\dbF^a$,
then in general the consistency condition \reff{consistent}
does not hold. For example, let $d=1$, $a_t := 1$, and
$b_t := 1+\1_{[1,\infty)}(t)$, $t\ge 0$,
as in Remark \ref{rem-augmentation0}.
Set $W^{\dbP^a}_t(\o):=B_t(\o) + \1_{(\O^a_1)^c}(\o)$ and
$W^{\dbP^b}_t(\o) := B_t(\o) + [B_t(\o)-B_{1}(\o)]\1_{[1,\infty)}(t)$.
Then both $W^{\dbP^a}$ and $W^{\dbP^b}$
are $\overline\dbF^{a}\cap \overline\dbF^{b}$-progressively measurable.
However, $\th^{a,b} = {1}$, but $\dbP^b(W^{\dbP^a}_0 = W^{\dbP^b}_0)
= \dbP^b(\O^a_1) = 0$, so we do not have
$W^{\dbP^a} = W^{\dbP^b}$, $\dbP^b$-almost surely on $[0,1]$.
\ep}
\end{rem}

\section{Quasi-sure stochastic analysis}
\label{sect-QS-Analysis}
\setcounter{equation}{0}

In this section, we fix again a separable class $\cA$ of diffusion coefficients generated by $\cA_0$, and set $\cP:=\{\dbP^a: a\in\cA\}$. We shall develop the $\cP$-quasi sure
stochastic analysis. We emphasize again that, when a probability measure $\dbP\in\cP$ is fixed, by Lemma \ref{lem-version} there is no need to distinguish the filtrations $\dbF^+, \dbF^\dbP$, and $\overline \dbF^\dbP$. 

We first introduce several spaces.
Denote by $\dbL^{0}$ the collection of all $\hat\cF^\cP_\infty$-measurable
random variables with appropriate dimension.
For each $p\in[1,\infty]$ and $\dbP\in\cP$,
we denote by $\dbL^{p}(\dbP)$ the corresponding
$\dbL^p$ space under the measure $\dbP$ and
\beaa
\hat\dbL^p &:=& \bigcap_{\dbP\in\cP} \dbL^{p}(\dbP).
\eeaa
Similarly, $\dbH^0:=\dbH^0(\dbR^d)$ denotes the collection
of all $\dbR^d$ valued $\hat\dbF^\cP$-progressively measurable processes.
$\dbH^{p}(\dbP^a)$ is the subset of all $H\in \dbH^0$
satisfying
\beaa
\|H\|_{T,\dbH^p(\dbP^a)}^p:=\dbE^{\dbP^a}\Big[\Big(\int_0^T
|a^{1\slash 2}_s H_s|^2 ds\Big)^{p/2}\Big]<\infty  &\mbox{for all}& T>0,
\eeaa
and $\dbH^{2}_{loc}(\dbP^a)$ is
the subset of $\dbH^0$ whose elements
satisfy $\int_0^T |a^{1\slash 2}_s H_s|^2  ds<\infty$,
 $\dbP^a$-almost surely, for all $T\ge 0$. Finally, we define
\beaa
\hat\dbH^p \;:=\; \bigcap_{\dbP\in\cP}\dbH^{p}(\dbP)
&\mbox{and}&
\hat\dbH^2_{loc} \;:=\; \bigcap_{\dbP\in\cP} \dbH^{2}_{loc}(\dbP).
\eeaa

The following two results are direct applications of
Theorem \ref{thm-aggregation}.
Similar results were also proved in \cite{DM,DHP}, see e.g.
 Theorem 2.1 in \cite{DM}, Theorem 36 in \cite{DHP}
and the Kolmogorov criterion of Theorem 31 in \cite{DHP}.

\begin{prop}[Completeness]\label{prop-Cauchy}
Fix $p\ge 1$, and let $\cA$ be satisfying \reff{Aseparable}.

\no
{\rm (i)}~Let $(X_n)_n\subset\hat\dbL^p$
be a Cauchy sequence under each $\dbP^a$, $a\in\cA$.
Then there exists a unique random variable
$X\in\hat\dbL^p$ such that $X_n\rightarrow X$ in
$\dbL^p(\dbP^a,\hat\cF^\cP_\infty)$ for every $a\in\cA$.

\no{\rm (ii)}~Let $(X_n)_n\subset\hat\dbH^p$ be a Cauchy
sequence under the norm $\|\cdot\|_{T,\dbH^p(\dbP^a)}$
for all $T\ge 0$ and $a\in\cA$. Then there exists a  unique process
$X\in \hat\dbH^p$ such that
$X_n\rightarrow X$ under the norm
$\|\cdot\|_{T,\dbH^p(\dbP^a)}$ for all $T\ge 0$ and $a\in\cA$.
\end{prop}

\proof
(i) By the completeness of
$\dbL^p(\dbP^a,\hat\cF^\cP_\infty)$,
we may find $X^a\in\dbL^p(\dbP^a,\hat\cF^\cP_\infty)$
such that $X_n\rightarrow X^a$ in
$\dbL^p(\dbP^a,\hat\cF^\cP_\infty)$.
The consistency condition of Theorem
\ref{thm-aggregation} is obviously satisfied
by the family $\{X^a,a\in\cA\}$, and the result follows.
 (ii) can be proved by a similar argument.
\ep

\begin{prop}[Kolmogorov continuity criteria]
\label{prop-Kolmogorov}
Let $\cA$ be satisfying \reff{Aseparable}, and $X$ be an $\hat\dbF^\cP$-progressively measurable process with values in $\dbR^n$.
We further assume that for some $p>1$, $X_t\in\hat\dbL^p$ for all $t\ge 0$ and satisfy
\beaa
\dbE^{\dbP^a}\left[|X_t-X_s|^p\right]
\le
c_a|t-s|^{n+\eps_a}
&\mbox{for some constants}&
c_a,\eps_a>0.
\eeaa
Then $X$ admits an $\hat\dbF^\cP$-progressively
measurable version $\tilde X$ which is
H\"older continuous, $\cP$-q.s.~{\rm(}with
H\"older exponent $\alpha_a<\eps_a/p$,
$\dbP^a$-almost surely for every $a\in\cA${\rm)}.
\end{prop}

\proof
We apply the Kolmogorov continuity criterion
under each $\dbP^a$, $a\in\cA$. This yields
a family of $\overline\dbF^{\dbP^a}$-progressively
measurable processes $\{X^a,a\in\cA\}$
such that $X^a=X$, $\dbP^a$-almost surely, and $X^a$ is H\"older continuous
with H\"older exponent $\alpha_a<\eps_a/p$,
$\dbP^a$-almost surely for every $a\in\cA$.
Also in view of Lemma \ref{lem-version},
we may assume without loss of generality
that $X^a$ is $\hat\dbF^\cP$-progressively measurable
for every $a\in\cA$. Since each $X^a$ is a
$\dbP^a$-modification
of $X$ for every $a\in\cA$, the consistency condition
of Theorem \ref{thm-aggregation} is immediately
satisfied by the family $\{X^a,a\in\cA\}$.
Then, the aggregated process $\tilde X$
constructed in that theorem has the
desired properties.\ep

\begin{rem}{\rm
The statements of Propositions \ref{prop-Cauchy}
and \ref{prop-Kolmogorov} can be weakened
further by allowing $p$ to depend on $a$.
\ep}
\end{rem}

We next construct the stochastic integral with
respect to the canonical process $B$ which is
simultaneously defined under all the mutually
singular measures $\dbP^a$, $a\in\cA$.  Such constructions
have been given in the literature but under
regularity assumptions on the integrand.
Here we only place standard conditions on
the integrand but not regularity.

\begin{thm}[Stochastic integration]
\label{thm-stochintegral}
For $\cA$ satisfying \reff{Aseparable}, let $H\in \hat\dbH^2_{loc}$ be given. Then, there exists a unique {\rm{(}}$\cP$-q.s.{\rm{)}} $\hat\dbF^\cP$-progressively measurable process $M$ such that $M$ is a local martingale under each $\dbP^a$ and
\beaa
M_t = \int_0^t H_s dB_s,~t\ge 0,
&\dbP^a\mbox{-almost surely for all}&
a\in\cA.
\eeaa
If in addition $H\in\hat\dbH^2$,
then for every $a\in\cA$,
$M$ is a square integrable
$\dbP^a$-martingale.  Moreover,
$\dbE^{\dbP^a}[M_t^2]=\dbE^{\dbP^a}
[\int_0^t |a^{1\slash 2}_s H_s|^2  ds]$ for all $t\ge 0$.
\end{thm}

\proof
For every $a\in\cA$, the stochastic integral
$M^a_t:=\int_0^t H_sdB_s$ is well-defined
$\dbP^a$-almost surely as a
$\overline\dbF^{\dbP^a}$-progressively measurable
process.  By Lemma \ref{lem-version}, we may
assume without loss of generality that
$M^a$ is $\hat\dbF^\cP$-adapted.
Following the arguments in Corollary \ref{cor-quadratic}, in particular by applying Lemma \ref{lem-Pconsistent}, it is clear
that the consistency condition \reff{consistent}
of Theorem \ref{thm-aggregation}
is  satisfied by the family $\{M^a,a\in\cA\}$.
Hence, there exists an aggregating process $M$.
The remaining statements in the theorem follows
from classical results for standard
stochastic integration under each $\dbP^a$.
\ep

We next study the martingale representation.

\begin{thm}[Martingale representation]
\label{thm-mgrep}
Let $\cA$ be a separable class of diffusion coefficients generated by
$\cA_0 \subset \cA_{\mrp}$. Let $M$ be an $\hat\dbF^\cP$-progressively measurable process which is a $\cP-$quasi sure local martingale,
that is, $M$ is a local martingale under $\dbP$ for all $\dbP\in\cP$. Then there exists a unique {\rm{(}}$\cP$-q.s.{\rm{)}} process $H\in \hat\dbH^2_{loc}$ such that
\beaa
M_t = M_0 +\int_0^t H_s dB_s, &t\ge 0,& \cP-\mbox{q.s.}.
\eeaa
\end{thm}

\proof By Proposition \ref{prop-cA},  $\cA\subset \cA_{\mrp}$.
Then for each $\dbP\in \cP$, all $\dbP-$martingales
can be represented as stochastic integrals
with respect to the canonical process.
Hence, there exists unique ($\dbP-$almost surely)
process $H^\dbP\in\dbH^{2}_{loc}(\dbP)$ such that
\beaa
M_t = M_0 +\int_0^t H^\dbP_s dB_s, &t\ge 0,&
\dbP\mbox{-almost surely.}
\eeaa
Then the quadratic covariation under $\dbP^b$ satisfies
\bea
\label{MB}
\la M, B\ra^{\dbP^b}_t = \int_0^t H^\dbP_s \hat a_s ds,\q &t\ge 0,& \dbP-\mbox{almost surely}.
\eea
Now for any $a, b\in \cA$, from the construction of quadratic covariation and  that of Lebesgue integrals, following similar arguments as in Corollary  \ref{cor-quadratic} one can easily check that
\beaa
\int_0^t H^{\dbP^a}_s \hat a_s ds = \la M, B\ra^{\dbP^a}_t=\la M, B\ra^{\dbP^b}_t =  \int_0^t H^{\dbP^b}_s \hat a_s ds, \q  0\le t<\th^{a,b},\q \dbP^a-\mbox{almost surely}.
\eeaa
This implies that
\beaa
H^{\dbP^a}\1_{[0, \th^{a,b})} =  H^{\dbP^b}\1_{[0, \th^{a,b})}, ~~ dt\times d\dbP^a-\mbox{almost surely}.
\eeaa
That is, the family $\{H^\dbP, \dbP\in\cP\}$ satisfies the
consistency condition \reff{consistent}.
Therefore,  we may aggregate them into a process $H$.
Then one may directly check that $H$ satisfies all the requirements.
\ep

\vspace{5mm}

There is also $\cP$-quasi sure decomposition of super-martingales.

\begin{prop}[Doob-Meyer decomposition]
\label{prop-decomposition}
For $\cA$ satisfying \reff{Aseparable}, assume an $\hat\dbF^\cP$-progressively measurable process $X$ is a $\cP$-quasi sure supermartingale, i.e.,
$X$ is a $\dbP^a$-supermartingale for all $a\in\cA$.
Then there exist a unique {\rm{(}}$\cP$-q.s.{\rm{)}}
$\hat\dbF^\cP$-progressively
measurable processes $M$ and $K$ such that $M$ is a
$\cP$-quasi sure local martingale and $K$ is
predictable and increasing, $\cP$-q.s., with
$M_0=K_0=0$, and $X_t = X_0 + M_t - K_t$, $t\ge 0$,
$\cP$-quasi surely.

If further $X$ is in class {\rm{(}}D{\rm{)}}, $\cP$-quasi surely,
i.e. the family $\{X_{\hat\t}, \hat\t\in\hat\cT\}$
is $\dbP$-uniformly integrable,
for all $\dbP\in\cP$, then $M$ is a $\cP$-quasi surely
uniformly integrable martingale.
\end{prop}

\proof For every $\dbP\in\cA$, we
apply Doob-Meyer decomposition theorem
(see e.g. Dellacherie-Meyer \cite{Dellacherie-Meyer}
Theorem VII-12).  Hence there exist a $\dbP$-local
martingale $M^\dbP$ and a $\dbP$-almost surely
increasing process $K^\dbP$ such that
$M^\dbP_0=K^\dbP_0=0$, $\dbP$-almost surely.
The consistency condition of Theorem \ref{thm-aggregation}
follows from the uniqueness of the Doob-Meyer decomposition.
Then, the aggregated processes provide the
universal decomposition.
\ep

\vspace{5mm}

The following results also follow from similar applications
of our main result.

\begin{prop}[It\^o's formula]\label{prop-Ito}
For $\cA$ satisfying \reff{Aseparable}, let $A, H$ be $\hat\dbF^\cP$-progressively measurable processes with values in
$\dbR$ and $\dbR^d$, respectively. Assume that $A$ has finite variation over each time interval $[0,t]$ and $H\in \hat\dbH^2_{loc}$. For $t \ge 0$, set $X_t:=A_t+\int_0^t H_sdB_s$. Then for any $C^2$ function $f:\dbR\rightarrow\dbR$,
we have
\beaa
f(X_t)
&=&
f(A_0)+\int_0^t f'(X_s)(dA_s+H_sdB_s)
+\frac12\int_0^t H_s^{\rm T}\hat a_s H_s f''(X_s)ds,~t\ge 0,
~\cP\mbox{-q.s..}
\eeaa
\end{prop}

\proof Apply It\^o's formula under each
$\dbP\in\cP$, and proceed as in the proof of
Theorem \ref{thm-stochintegral}.
\ep

\begin{prop}[local time]\label{cor-localtime}
For $\cA$ satisfying \reff{Aseparable}, let $A$, $H$ and $X$ be as in
Proposition \ref{prop-Ito}.
Then for any $x\in\dbR$, the local time $\{L^x_t,t\ge 0\}$ exists
$\cP$-quasi surely and is given by,
\beaa
2L_t^x
&=&
|X_t-x|-|X_0-x|-\int_0^t {\rm sgn}
(X_s-x)(dA_s+H_sdB_s),~t\ge 0,~\cP-\mbox{q.s..}
\eeaa
\end{prop}

\proof
Apply Tanaka's formula under each $\dbP\in\cP$
and proceed as in the proof of Theorem \ref{thm-stochintegral}.
\ep

\vspace{5mm}

Following exactly as in the previous results, we
obtain a Girsanov theorem in this context as well.

\begin{prop}[Girsanov]
For $\cA$ satisfying \reff{Aseparable}, let $\phi$ be $\hat\dbF^\cP$-progressively measurable and $\int_0^t |\phi_s|^2 ds <\infty$ for all $t\ge 0$, $\cP$-quasi surely. Let
\beaa
Z_t:=\exp{\left(\int_0^t\phi_sdW_s-\frac12\int_0^t|\phi_s|^2ds\right)}
&\mbox{and}&
\tilde{W}_t:=W_t-\int_0^t\phi_sds,~t\ge 0,
\eeaa
where $W$ is the $\cP$-Brownian motion of \reff{W}.
Suppose that for each $\dbP \in \cP$,
$\dbE^{\dbP}[Z_T]=1$ for some $T\ge 0$.
On $\hat\cF_T$ we define the probability measure
$\dbQ^\dbP$ by $d\dbQ^\dbP=Z_Td\dbP$. Then,
\beaa
\dbQ^\dbP\circ \tilde W^{-1} = \dbP\circ W^{-1}
&\mbox{for every}&
\dbP\in\cP,
\eeaa
i.e. $\tilde W$ is a $\dbQ^\dbP$-Brownian
motion on $[0,T,]$ for every $\dbP\in\cP$.
\end{prop}
\vspace{5mm}

We finally discuss stochastic differential equations
in this framework. Set $\bQ^m:=\{(t,\bx):t\ge 0,\bx\in C[0,t]^m\}$.
Let $b,\si$ be two functions from $\Omega\times \bQ^m$ to $\dbR^m$ and
$\cM_{m,d}(\dbR)$, respectively.
Here, $\cM_{m,d}(\dbR)$ is the space of $m\times d$
matrices with real entries.
We are interested in the problem of solving
the following stochastic differential equation
{\it simultaneously under all} $\dbP\in\cP$,
\bea\label{sde}
X_t=X_0+\int_0^t b_s(\underline{X}_s)
ds+\int_0^t\sigma_s(\underline{X}_s)dB_s,
&t\ge 0,&
\cP-\mbox{q.s.,}
\eea
where $\underline{X}_t:=(X_s,s\le t)$.

\begin{prop}\label{prop-sde}
Let $\cA$ be satisfying \reff{Aseparable}, and assume that, for every
$\dbP\in\cP$ and  $\t\in\cT$, the equation \reff{sde} has a unique
$\dbF^{\dbP}$-progressively measurable strong solution on interval $[0,\t]$. Then there is a $\cP$-quasi surely aggregated  solution to \reff{sde}.
\end{prop}

\proof
For each $\dbP\in\cA$, there is a
$\dbP$-solution $X^\dbP$ on $[0,\infty)$,
which we may consider in its
$\hat\dbF^\cP$-progressively measurable
version by Lemma \ref{lem-version}. The uniqueness
on each $[0,\t]$,$\t\in\cT$ implies that the family
$\{X^\dbP,\dbP\in\cP\}$ satisfies the consistency
condition of Theorem \ref{thm-aggregation}.
\ep

\section{An application}
\setcounter{equation}{0}
\label{sec-UVM}

As an application of our theory,
we consider the problem of super-hedging
contingent claims under volatility uncertainty,
which was studied by Denis and Martini \cite{DM}.
In contrast with their approach, our framework allows
to  obtain the existence of the optimal hedging strategy.
However, this is achieved at the price of restricting
the non-dominated family of probability measures.

{ We also mention a related recent paper by Fernholz and Karatzas \cite{FK}
whose existence results are obtained in the Markov case with a continuity assumption on the corresponding value function.}

Let $\cA$ be a separable class of diffusion coefficients
generated by $\cA_0$, and $\cP:=\{\dbP^a:a\in\cA\}$
be the corresponding family of measures.
We consider a fixed time horizon, say $T=1$.
Clearly all the results in previous sections can be extended to this setting, after some obvious modifications.
Fix a nonnegative $\hat\cF_1-$measurable
real-valued random variable $\xi$.
The superhedging cost of $\xi$ is defined by
\beaa
v(\xi) &:=& \inf\left\{x:~x+\int_0^1H_sdB_s
\ge \xi,~~\cP\mbox{-q.s. for some}~~H\in\cH\right\},
\eeaa
where the stochastic integral $\int_0^\cd H_s dB_s$ is defined in the sense of Theorem \ref{thm-stochintegral} and $ H \in \dbH^0$ belongs to $\cH$ if and only if
\beaa
\int_0^1 H_t^T\hat a_t H_tdt<\infty~\cP\mbox{-q.s.}
                           &\mbox{and}&
\int_0^.H_sdB_s ~\mbox{is a $\cP$-q.s.~supermartingale}.
\eeaa

We shall provide a dual formulation of
the problem $v(\xi)$ in terms of the following
dynamic optimization problem,
\bea
\label{Vat}
V^{\dbP^a}_{\hat\t}
&:=&
\esup_{b\in\cA(\hat\t,a)}~^{\!\!\!\!\!\!\dbP^a}
       \dbE^{\dbP^{b}}[\xi|\hat\cF_{\hat\t}],~~\dbP^a\mbox{-a.s.},~~
       a\in\cA,~\hat\t\in\hat\cT,
\eea
where
$$\cA(\hat\t,a) := \{b\in \cA:
\th^{a,b} > \hat\t~~\mbox{or}~~\th^{a,b}= \hat\t=1\}.
$$

\begin{thm}\label{thm-DM}
Let $\cA$ be a separable class of diffusion coefficients generated by $\cA_0\subset \cA_{\mrp}$. Assume that the family of random variables
$\{V^\dbP_{\hat\t}, \hat\t\in\hat\cT\}$ is uniformly integrable under all $\dbP\in\cP$. Then
\be
\label{v=V}
v(\xi) = V(\xi) \;:=\; \sup_{a\in\cA} \|V^{\dbP^a}_0\|_{\dbL^\infty(\dbP^a)}.
\ee
Moreover, if $v(\xi)<\infty$, then
there exists $H\in\cH$ such that
$v(\xi)+\int_0^1 H_sdB_s\ge\xi$, $\cP$-q.s..
\end{thm}

\bs

To prove the theorem, we need the following (partial)
dynamic programming principle.

\begin{lem}
\label{lem-DPP}
Let $\cA$ be satisfying \reff{Aseparable}, and assume $V(\xi)<\infty$. Then, for any $\hat\t_1, \hat\t_2\in\hat\cT$ with $\hat\t_1\le \hat\t_2$,
$$
V^{\dbP^a}_{\hat\t_1} \ge \dbE^{\dbP^{b}}
\big[V^{\dbP^b}_{\hat\t_2}|\hat\cF_{\hat\t_1}\big],
\dbP^a\mbox{-almost surely for all}\
a\in \cA~\mbox{and}~b\in \cA(a,\hat\t_1).
$$
\end{lem}

\proof By the definition of essential supremum,
see e.g. Neveu \cite{Neveu} (Proposition VI-1-1), there exist a
sequence $\{b_j, j\ge 1\}\subset \cA(b, \hat\t_2)$
such that $V^{\dbP^b}_{\hat\t_2} =\sup_{j\ge 1}
\dbE^{\dbP^{b_j}}[\xi|\hat\cF_{\hat\t_2}]$,
$\dbP^{b}$-almost surely. For $n\ge 1$, denote
$V^{b,n}_{\hat\t_2} := \sup_{1\le j\le n}
\dbE^{\dbP^{b_j}}[\xi|\hat\cF_{\hat\t_2}]$.
Then $V^{b,n}_{\hat\t_2}\uparrow  V^{\dbP^b}_{\hat\t_2}$,
$\dbP^{b}$-almost surely as $n\to\infty$.
By the monotone convergence theorem,
we also have $\dbE^{\dbP^{b}}[V^{b,n}_{\hat\t_2}|
\hat\cF_{\hat\t_1}]\uparrow
\dbE^{\dbP^{b}}[V^{\dbP^b}_{\hat\t_2}|
\hat\cF_{\hat\t_1}]$, $\dbP^{b}$-almost surely, as $n\to\infty$.
Since $b\in\cA(a, \hat\t_1)$, $\dbP^{b} = \dbP^{a}$ on
$\hat\cF_{\hat\t_1}$. Then
$\dbE^{\dbP^{b}}[V^{b,n}_{\hat\t_2}|\hat\cF_{\hat\t_1}]
\uparrow \dbE^{\dbP^{b}}
[V^{\dbP^b}_{\hat\t_2}|\hat\cF_{\hat\t_1}]$,
$\dbP^{a}$-almost surely, as $n\to\infty$.
Thus it suffices to show that
\bea
\label{Vnn}
V^{\dbP^a}_{\hat\t_1} \ge
\dbE^{\dbP^{b}}[V^{b,n}_{\hat\t_2}|\hat\cF_{\hat\t_1}],
&& \dbP^{a}\mbox{-almost surely for all}~n\ge 1.
\eea

Fix $n$ and define
\beaa
\th^b_n := \min_{1\le j\le n} \th^{b, b_j}.
\eeaa
By Lemma \ref{lem-Pconsistent}, $\th^{b, b_j}$ are
$\dbF$-stopping times taking only countably
many values, then so is $\th^{b}_n$. Moreover,
since $b_j\in\cA(b, \hat\t_2)$, we have either
$\th^b_n>\hat\t_2$ or $\th^{b}_n=\hat\t_2=1$.
Following exactly the same arguments as in the
proof of \reff{inclusion}, we arrive at
\beaa
\hat\cF_{\hat\t_2} \subset \left(\cF_{\th^b_n}
\vee \cN^{\dbP^b}(\cF_1)\right).
\eeaa
Since $\dbP^{b_j} = \dbP^b$ on $\hat\cF_{\hat\t_2}$,
without loss of generality we may assume  the
random variables $\dbE^{\dbP^{b_j}}[\xi|\hat\cF_{\hat\t_2}]$
and  $V^{b,n}_{\hat\t_2}$ are $\cF_{\th^b_n}$-measurable.
Set $A_j := \{\dbE^{\dbP^{b_j}}[\xi|\hat\cF_{\hat\t_2}]
=V^{b,n}_{\hat\t_2}\}$ and
$\tilde A_1 := A_1$, $\tilde A_j :=
A_j\backslash \cup_{i<j} A_i$, $2\le j\le n$.
Then $\tilde A_1, \cds, \tilde A_n$ are
$\cF_{\th^b_n}$-measurable and form a
partition of $\O$.  Now set
$$
\tilde b(t) := b(t)\1_{[0,\hat\t_2)}(t) + \sum_{j=1}^n b_j(t)\1_{\tilde A_j}\1_{[\hat\t_2, 1]}(t).
$$
We claim that $\tilde b \in \cA$.  Equivalently, we need
to show that $\tilde b$ takes the form \reff{Aa}.
We know that $b$ and $b_j$ have the form
\beaa
b(t) = \sum_{m=0}^\infty \sum_{i=1}^\infty
b^{0,m}_i \1_{E^{0,m}_i} \1_{[\t^0_m, \t^0_{m+1})}
&\mbox{and}&
b_j(t) = \sum_{m=0}^\infty \sum_{i=1}^\infty
b^{j,m}_i \1_{E^{j,m}_i} \1_{[\t^j_m, \t^j_{m+1})}
\eeaa
with the stopping times and sets as before.
Since $b_j(t) = b(t)$ for $t\le \th^b_n$ and $j=1,\cds,n$,
\beaa
\tilde b(t) &=& b(t) \1_{[0, \th^{b}_n)} +
\sum_{j=1}^n \1_{\tilde A_j} b_j(t)\1_{[\th^{b}_n, 1]}(t)\\
&=& \sum_{m=0}^\infty \sum_{i=1}^\infty b^{0,m}_i
\1_{E^{0,m}_i\cap \{\t^0_m<\th^b_n\}}
\1_{[\t^0_m\wedge \th^b_n, \t^0_{m+1}\wedge\th^b_n)}\\
&&+ \sum_{j=1}^n\sum_{m=0}^\infty
\sum_{i=1}^\infty b^{j,m}_i \1_{E^{j,m}_i\cap
\tilde A_j\cap \{\t^j_{m+1}>\th^b_n\}}
 \1_{[\t^j_m\vee \th^b_n, \t^j_{m+1}\vee \th^b_n)}.
\eeaa
By Definition \ref{defn-cA}, it is clear that
$\t^0_{m} \wedge \th^b_n$ and $\t^j_m\vee \th^b_n$
are $\dbF$-stopping times and take only countably many values,
for all $m\ge 0$ and $1\le j\le n$.
For $m\ge 0$ and $1\le j\le n$, one can easily see that
$E^{0,m}_i\cap \{\t^0_m<\th^b_n\}$ is
$\cF_{\t^0_m\wedge \th^b_n}$-measurable and that
$E^{j,m}_i\cap \tilde A_j\cap \{\t^j_{m+1}>\th^b_n\}$ is
$\cF_{\t^j_m\vee \th^b_n}$-measurable.
By ordering the stopping times
$\t^0_{m} \wedge \th^b_n$ and $\t^j_m\vee \th^b_n$
we prove our claim that
$\tilde b \in \cA$.

It is now clear that
$\tilde b\in\cA(b, \hat\t_2)\subset\cA(a,\hat\t_1)$.
Thus,
\beaa
V^{\dbP^a}_{\hat\t_1}
&\ge& \dbE^{\dbP^{\tilde b}}[\xi|\hat\cF_{\hat\t_1}]
= \dbE^{\dbP^{\tilde b}}\Big[\dbE^{\dbP^{\tilde b}}
[\xi|\hat\cF_{\hat\t_2}]\Big|\hat\cF_{\hat\t_1}\Big]\\
&=& \dbE^{\dbP^{\tilde b}}\Big[\sum_{j=1}^n
\dbE^{\dbP^{\tilde b}}[\xi\1_{\tilde A_j}
|\hat\cF_{\hat\t_2}]\Big|\hat\cF_{\hat\t_1}\Big]\\
&=& \dbE^{\dbP^{\tilde b}}\Big[\sum_{j=1}^n
\dbE^{\dbP^{b_j}}[\xi\1_{\tilde A_j}
|\hat\cF_{\hat\t_2}]\Big|\hat\cF_{\hat\t_1}\Big]\\
&=& \dbE^{\dbP^{\tilde b}}\Big[\sum_{j=1}^n
V^{b,n}_{\hat\t_2} \1_{\tilde A_j}
\Big|\hat\cF_{\hat\t_1}\Big]
=\dbE^{\dbP^{\tilde b}}[V^{b,n}_{\hat\t_2}
\Big|\hat\cF_{\hat\t_1}],~~\dbP^{a}\mbox{-almost surely.}
\eeaa
Finally, since $\dbP^{\tilde b} = \dbP^b$ on
$\hat\cF_{\hat\t_2}$ and $\dbP^b=\dbP^a$ on
$\hat\cF_{\hat\t_1}$, we prove \reff{Vnn} and hence the lemma.
\ep

\bs

{\it Proof of Theorem \ref{thm-DM}.}
We first prove that $v(\xi)\ge V(\xi)$.
If $v(\xi)=\infty$, then the inequality is obvious.
If $v(\xi)<\infty$, there are $x\in\dbR$ and $H\in\cH$
such that the process $X_t:=x+\int_0^t H_sdB_s$
satisfies $X_1\ge\xi$, $\cP-$quasi surely.
Notice that the process $X$ is a
$\dbP^b$-supermartingale for every $b\in\cA$. Hence
$$
x=X_0\ge\dbE^{\dbP^b}[X_1|\hat\cF_0]
\ge\dbE^{\dbP^b}[\xi|\hat\cF_0], \q \dbP^b-\mbox{a.s.} \q \forall \ b\in\cA.
$$
By Lemma \ref{lem-Pconsistent}, we know that
$\dbP^a = \dbP^b$ on $\hat\cF_0$ whenever
$a\in \cA$ and $b\in \cA(0,a)$.  Therefore,
\beaa
x\ge\dbE^{\dbP^b}[\xi|\hat\cF_0],
~~\dbP^a\mbox{-a.s..}
\eeaa
The definition of $V^{\dbP^a}$ and the above inequality
imply that
$x\ge V^{\dbP^a}_0$, $\dbP^a$-almost surely.
This implies that $x\ge \|V^{\dbP^a}_0\|_{\dbL^\infty(\dbP^a)}$
for all $a\in\cA$.  Therefore, $x\ge V(\xi)$. Since this holds
for any initial data $x$ that is super-replicating $\xi$,
we conclude that $v(\xi)\ge V(\xi)$.

We next prove the opposite inequality.
Again, we may assume that $V(\xi)<\infty$.
Then $\xi\in\hat\dbL^1$.
For each $\dbP\in\cP$, by Lemma \ref{lem-DPP} the family $\{V^\dbP_{\hat\t}, \hat\t\in\hat\cT\}$ 
satisfies the (partial) dynamic programming principle. Then following standard arguments 
(see e.g. \cite{EQ} Appendix A2), we construct from this family a {\cad}
$(\hat\dbF^\cP,\dbP)$-supermartingale $\hat V^\dbP$ defined by,
\bea
\label{Va}
\hat V^\dbP_t := \limsup_{\dbQ\ni r\downarrow t}
V^\dbP_r,~~t\in [0,1].
\eea

Also for each $\hat\t\in\hat\cT$, it is clear that the family
$\{V^\dbP_{\hat\t},\dbP\in\cP\}$ satisfies the consistency condition
\reff{consistent1}. Then it follows immediately from
\reff{Va} that $\{\hat V^\dbP_t, \dbP\in\cP\}$ satisfies
the consistency condition \reff{consistent1} for all $t\in [0,1]$.
Since $\dbP$-almost surely $\hat V^\dbP$ is {\cad},
the family of processes $\{\hat V^\dbP, \dbP\in\cP\}$
also satisfy the consistency condition \reff{consistent}.
We then conclude from Theorem \ref{thm-aggregation}
that there exists a unique aggregating process $\hat V$.

Note that $\hat V$ is a $\cP$-quasi sure supermartingale.
Then it follows from the Doob-Meyer decomposition of
Proposition \ref{prop-decomposition} that
there exist a $\cP$-quasi sure local martingale $M$
and a $\cP$-quasi sure increasing process $K$
such that $M_0=K_0=0$ and $\hat V_t = \hat V_0
+ M_t - K_t$, $t\in [0,1)$, $\cP$-quasi surely.
Using the uniform integrability hypothesis of this theorem,
we conclude that the previous decomposition
holds on $[0,1]$ and the process $M$ is a $\cP$-quasi
sure martingale on $[0,1]$.

In view of the martingale representation
Theorem \ref{thm-mgrep}, there exists an
$\hat\dbF^\cP$-progressively measurable process
$H$ such that $\int_0^1 H_t^T\hat a_t H_tdt<\infty$
and $\hat V_t = \hat V_0 + \int_0^t H_s dB_s - K_t$, $t\ge 0$,
$\cP$-quasi surely. Notice that $\hat V_1 = \xi$ and
$K_1\ge K_0=0$.  Hence  $ \hat V_0 + \int_0^1 H_s dB_s \ge \xi$,
$\cP$-quasi surely. Moreover, by the definition of
$V(\xi)$, it is clear that $V(\xi) \ge \hat V_0$, $\cP$-quasi surely.
Thus $V(\xi) + \int_0^1 H_s dB_s \ge \xi$, $\cP$-quasi surely.

Finally, since $\xi$ is nonnegative, $\hat V \ge 0$.
Therefore,
$$
V(\xi) + \int_0^t H_s dB_s
\ge \hat V_0 + \int_0^t H_s dB_s
\ge \hat V_t \ge 0, \qquad \cP-q.s..
$$
This implies that $H\in\cH$, and thus $V(\xi)\ge v(\xi)$.
\ep

\vspace{5mm}

\begin{rem}
\label{rem-abound}{\rm
Denis and Martini \cite{DM} require
\bea
\label{cA0bound}
\underline{a}\le a\le \overline{a}
&\mbox{for all}& a\in\cA,
\eea
for some given constant matrices
$\underline{a}\le\overline{a}$ in $\dbS^{>0}_d$.
We do not impose this constraint.
In other words, we may allow
$\underline a = 0$ and $\overline a = \infty$.
Such a relaxation is important in problems of static
hedging in finance, see e.g.~\cite{CL} and the references therein.
However, we still require that each $a\in \cA$ takes values in $\dbS^{>0}_d$.
\ep}\end{rem}

We shall introduce the set $\cA_S\subset \cA_{\mrp}$ induced from
strong formulation in Section \ref{ss.cPS}. When $\cA_0\subset \cA_S$,
we have the following additional interesting properties.

\begin{rem}
\label{rem-Blumenthal}{\rm
If each $\dbP\in\cP$ satisfies the Blumenthal zero-one law
(e.g.~if $\cA_0\subset \cA_S$ by Lemma \ref{lem-mgarep} below), then $V^{\dbP^a}_0$ is a constant
for all $a\in\cA$, and thus \reff{v=V} becomes
\beaa
v(\xi) = V(\xi) := \sup_{a\in\cA} V^{\dbP^a}_0.
\eeaa}
\end{rem}
\begin{rem}
\label{rem-dense}{\rm
In general, the value $V(\xi)$ depends on $\cA$,
then so does $v(\xi)$. However, when $\xi$ is uniformly
continuous in $\o$ under the uniform norm, we show in \cite{STZ09c} that
\bea
\label{duality-cPS}
\sup_{\dbP\in\overline{\cP}_S} \dbE^{\dbP}[\xi] =
\inf\left\{x:~x+\int_0^1H_sdB_s
\ge \xi,~~\dbP\mbox{-a.s. for all}~\dbP\in\overline\cP_S,
~~\mbox{for some}~H\in\overline\cH\right\},
\eea
and the optimal superhedging strategy $H$ exists,
where $\overline\cH$ is the space of $\dbF$-progressively
measurable $H$ such that, for all $\dbP\in\overline\cP_S$,
$\int_0^1 H_t^T\hat a_t H_tdt<\infty$, $\dbP$-almost surely
and $\int_0^.H_sdB_s$ is a $\dbP$-supermartingale.
Moreover, if $\cA\subset \cA_S$ is dense in some sense, then
\beaa
V(\xi) = v(\xi) = \mbox{the $\overline\cP_S$-superhedging
cost in \reff{duality-cPS}}.
\eeaa
In particular, all functions are independent of the
choice of $\cA$.
This issue is discussed in details in
our accompanying paper \cite{STZ09c} (Theorem 5.3 and Proposition 5.4), 
where we establish a duality result for a
more general setting called the second order target problem.
However, the set-up in \cite{STZ09c} is more general
and this independence can be proved by the above
arguments under suitable assumptions. 
\ep}
\end{rem}

\section{Mutually singular measures induced by strong formulation }
\label{ss.cPS}

We recall the set $\overline{\cP}_S$ introduced in the Introduction as
\bea
\label{overlinecPS}
\overline{\cP}_S :=\left\{\dbP^\a_S: \a \in \overline{\cA}\right\} &\mbox{where}&
\dbP^\a_S := \dbP_0 \circ \left(X^\a\right)^{-1},
\eea
and $X^\a$ is given in \reff{e.xa}.
Clearly $\overline\cP_S\subset \overline{\cP}_W$.
Although we do not use it in the present paper, this class is
important both in theory and in applications.
We remark that Denis-Martini \cite{DM} and
our paper \cite{STZ09b} consider
the class $\overline{\cP}_W$ while Denis-Hu-Peng \cite{DHP}
and our paper \cite{STZ09d}
consider the class $\overline{\cP}_S$,
up to some technical restriction of the
diffusion coefficients.

We start the analysis of this set by noting that
\bea
\label{BXa}
\a ~\mbox{is the  quadratic variation density of}~X^\a
~~\mbox{and}~~dB_s= \a^{{-1\slash 2}}_s dX^\a_s, \ \
{\mbox{under}}\ \dbP_0.
\eea
Since $B$ under $\dbP^\a_S$ has the same
distribution as $X^\a$ under $\dbP_0$, it is clear that
\bea
\label{e.equivalent}
\mbox{the}~\dbP^\a_S{\mbox{-distribution of}}\
(B, \hat a, W^{\dbP^\a_S})\
{\mbox{is equal to the}}\
\dbP_0{\mbox{-distribution of}}\ (X^\a, \a, B).
\eea
In particular, this implies that
\bea
\label{strong}
\ba{c}
\hat a(X^\a) = \a(B),~~\dbP_0\mbox{-a.s.},\q
\hat a(B) = \a(W^{\dbP^\a_S}),~~\dbP^\a_S\mbox{-a.s.},\\
\mbox{and for any}~ a\in \overline\cA_W(\dbP^\a_S),
~ X^\a ~\mbox{is a strong solution to SDE \reff{SDE0} with coefficient}~a.
\ea
\eea
Moreover we have the following characterization of
$\overline{\cP}_S$ in terms of the filtrations.
\begin{lem}
\label{l.technical}
$\dis
\overline{\cP}_S = \left\{\dbP\in\overline{\cP}_W:
\overline{\dbF^{W^\dbP}}^\dbP = \overline{\dbF}^\dbP\right\}.
$
\end{lem}
\proof By \reff{BXa}, $\a$ and $B$ are
$\overline{\dbF^{X^\a}}^{\dbP_0}$-progressively measurable.
Since $\dbF$ is generated by $B$, we conclude that
$\dbF \subset \overline{\dbF^{X^\a}}^{\dbP_0} $.
By completing the filtration we next obtain that
$ \overline{\dbF}^{\dbP_0} \subset \overline{\dbF^{X^\a}}^{\dbP_0} $.
Moreover, for any $\a \in \overline{\cA}$,
it is clear that $\dbF^{X^\a}\subset \overline{\dbF}^{\dbP_0}$.
Thus, $\overline{\dbF^{X^\a}}^{\dbP_0} = \overline{\dbF}^{\dbP_0}$.
Now, we invoke \reff{e.equivalent} and conclude
$\overline{\dbF^{W^\dbP}}^\dbP = \overline{\dbF}^\dbP$
for any $\dbP=\dbP^\a_S\in\overline{\cP}_S $.

Conversely, suppose $\dbP\in\overline{\cP}_W$ be such that
$\overline{\dbF^{W^\dbP}}^\dbP = \overline{\dbF}^\dbP$.
Then $B = \b(W^\dbP_\cd)$ for some measurable mapping
$\b: \bQ\to \dbS^{>0}_d$. Set $\a := \b(B_\cd)$,
we conclude  that $\dbP = \dbP^\a_S$.
\ep

\vspace{5mm}

The following result shows that the measures
$\dbP\in\overline{\cP}_S$ satisfy MRP and the  Blumental zero-one law.

\begin{lem}\label{lem-mgarep}
$\overline{\cP}_S\subset \overline{\cP}_{\mrp}$
and every $\dbP\in \overline{\cP}_S$ satisfies the
Blumenthal zero-one law.
\end{lem}

\proof Fix $\dbP\in \overline{\cP}_S$.  We first show
that $\dbP\in \overline{\cP}_{\mrp}$.   Indeed, for any
$(\overline\dbF^{\dbP},\dbP)$-local martingale $M$,
Lemma \ref{l.technical} implies that $M$ is a
$(\overline{\dbF^{W^\dbP}}^{\dbP},\dbP)$-local martingale.
Recall that $W^\dbP$ is a $\dbP$ Brownian motion.
Hence, we now can use the standard martingale representation theorem.
Therefore,  there exists a unique
$\overline{\dbF^{W^\dbP}}^{\dbP}$-progressively
measurable process $\tilde H$ such that
\beaa
\int_0^t |\tilde H_s|^2 ds<\infty ~~\mbox{and}~~
M_t = M_0+\int_0^t \tilde H_s d W^\dbP_s, ~~t\ge 0,
~~\dbP\mbox{-a.s..}
\eeaa
Since $\hat a>0$, $dW^\dbP=  \hat a^{-{1\slash 2}} dB$.
So one can check directly that the process
$H := \hat a^{-{1\slash 2}}\tilde H$ satisfies all the requirements.

We next prove the Blumenthal zero-one law. For this purpose
fix  $E\in \cF_{0+}$.  By Lemma \ref{l.technical},
$E\in \overline{\cF^{W^\dbP}_0}^{\dbP}$.
Again we recall that $W^\dbP$ is a $\dbP$ Brownian motion
and use the standard Blumenthal zero-one law for the Brownian motion.
Hence $\dbP(E)\in \{0,1\}$.
\ep

\vspace{5mm}

We now define analogously the following spaces of measures and diffusion processes.
\bea
\label{cPS}
\cP_S := \overline\cP_S \cap \cP_W, && \cA_S := \left\{a\in \cA_W: \dbP^a \in \cP_S\right\}.
\eea
Then it is clear that
\beaa
\cP_S \subset \cP_{\mrp}\subset \cP_W &\mbox{and}& \cA_S\subset \cA_{\mrp}  \subset \cA_W.
\eeaa
The conclusion $\cP_S \subset \cP_W$ is strict, see Barlow \cite{Barlow}. We remark that one can easily check that the diffusion process $a$ in Examples \ref{eg-lip} and \ref{eg-Piecewise constant} and the generating class $\cA_0$ in Examples \ref{eg-deterministic}, \ref{eg-Lipschitz}, and \ref{eg-hP} are all in $\cA_S$.

Our final result extends Proposition \ref{prop-cA}.

\begin{prop}
\label{prop-cA2}
Let $\cA$ be a separable class of
diffusion coefficients generated by
$\cA_0$. If $\cA_0\subset \cA_S$, then $\cA\subset\cA_S$.
\end{prop}

\proof Let $a$ be given in
the form \reff{Aa} and,  by Proposition \ref{prop-cA}, $\dbP$ be the unique weak solution to SDE \reff{SDE0} on $[0,\infty)$ with coefficient $a$ and initial condition $\dbP(B_0=0)=1$. 
By Lemma \ref{l.technical} and its proof, it suffices to show that  $a$ is  $\overline {\dbF^{W^\dbP}}^\dbP$-adapted.  Recall \reff{Aa}. We prove by induction on $n$ that
\bea
\label{a-adapt}
a_t \1_{\{t<\t_n\}}~~ \mbox{is}~~ \overline {\cF^{W^\dbP}_{t\wedge\t_n}}^\dbP-\mbox{measurable for all}~~t\ge 0.
\eea

Since $\t_0=0$, $a_0$ is $\cF_0$-measurable,  and $\dbP(B_0=0)=1$, \reff{a-adapt} holds when $n=0$.  Assume  \reff{a-adapt}  holds true for $n$. Now we consider $n+1$. Note that
\beaa
a_t  \1_{\{ t<\t_{n+1}\}} = a_t  \1_{\{ t<\t_{n}\}}  +a_t  \1_{\{\t_n\le t<\t_{n+1}\}} .
\eeaa
By the induction assumption it suffices to show that
\bea
\label{a-adapt1}
a_t \1_{\{ t<\t_{n+1}\}}~~ \mbox{is}~~ \overline {\cF^{W^\dbP}_{\t_n \vee t\wedge\t_{n+1}}}^\dbP-\mbox{measurable for all}~~t\ge 0.
\eea
Apply Lemma \ref{lem-Astructure}, we have $a_t = \sum_{m\ge 1} a_m(t)\1_{E_m}$ for $t<\t_{n+1}$, where $a_m\in \cA_0$ and $\{E_m, m\ge 1\}\subset \cF_{\t_n}$ form a partition of $\O$.  Let $\dbP^m$ denote the unique weak solution to SDE \reff{SDE0} on $[0,\infty)$ with coefficient $a_m$ and initial condition $\dbP^m(B_0=0)=1$. Then by Lemma \ref{lem-Pconsistent} we have, for each $m\ge 1$, 
\bea
\label{P=Pm}
\dbP(E\cap E_m) = \dbP^m(E\cap E_m),\q \forall E\in \cF_{\t_{n+1}}.
\eea
Morover, by (\ref{WP}) it is clear that 
\bea
\label{WP=WPm}
W^\dbP_t = W^{\dbP^m}_t, ~~0\le t\le \t_{n+1}~~, \dbP-\mbox{a.s. on}~~E_m ~~(\mbox{and} ~~\dbP^m-\mbox{a.s. on}~~E_m).
\eea
Now since $a_m\in \cA_0 \subset \cA_S$, we know $a_m(t)\1_{\{t<\t_{n+1}\}}$ is $\overline {\cF^{W^{\dbP^m}}_{t\wedge\t_{n+1}}}^{\dbP^m}-$measurable. This, together with the fact that $E_m\in \cF_{\t_n}$, implies that $a_m(t)\1_{\{ t<\t_{n+1}\}}\1_{E_m}$ is  $\overline {\cF^{W^{\dbP^m}}_{\t_n\vee t\wedge\t_{n+1}}}^{\dbP^m}-$measurable.
By (\ref{P=Pm}), (\ref{WP=WPm}) and that $a_t = a_m(t)$ for $t<\t_{n+1}$ on $E_m$,  we see that $a_t\1_{\{ t<\t_{n+1}\}}\1_{E_m}$ is  $\overline {\cF^{W^{\dbP}}_{\t_n\vee t\wedge\t_{n+1}}}^{\dbP}-$measurable. Since $m$ is arbitrary, we get
\beaa
a_t\1_{\{ t<\t_{n+1}\}} = \sum_{m\ge 1}a_t\1_{\{ t<\t_{n+1}\}}\1_{E_m}
\eeaa
is $\overline {\cF^{W^{\dbP}}_{\t_n\vee t\wedge\t_{n+1}}}^{\dbP}-$measurable. This proves \reff{a-adapt1}, and hence the proposition.
\ep

\section{Appendix}
\setcounter{equation}{0}

In this Appendix we provide a few more
examples concerning weak solutions of
\reff{SDE0} and complete the remaining technical proofs.

\subsection{Examples}

\begin{eg} {\rm (No weak solution)
Let $a_0 = 1$, and for $t>0$,
\beaa
a_t := 1 + \1_E,&\mbox{where}& E:=
\Big\{\limsup_{h\downarrow 0}{B_{h}-B_0\over \sqrt{2h\ln\ln h^{-1}}} \neq 2\Big\}.
\eeaa
Then $E\in\cF_{0+}$. Assume $\dbP$ is
a weak solution to (\ref{SDE0}). On $E$, $a=2$, then
$\limsup_{h\downarrow 0}{B_{h}-B_0\over \sqrt{2h\ln\ln h^{-1}}}
= 2$, $\dbP$-almost surely, thus $\dbP(E)=0$.
On $E^c$, $a=1$, then $\limsup_{h\downarrow 0}{B_{h}-
B_0\over \sqrt{2h\ln\ln h^{-1}}} = 1$, $\dbP$-almost surely and
thus $\dbP(E^c)=0$. Hence there can not be any weak solutions.
\ep
}\end{eg}

\begin{eg} {\rm(Martingale measure without Blumenthal 0-1 law)
Let $\O':= \{1,2\}$ and $\dbP_0'(1) = \dbP_0'(2) = {1\over 2}$.
Let $\tilde\O:= \O\times\O'$ and $\tilde\dbP_0$ the
product of $\dbP_0$ and $\dbP_0'$. Define
\beaa
 \tilde B_t(\o,1) := \o_t,\q  \tilde B_t(\o,2) := 2\o_t.
\eeaa
Then $\tilde\dbP:= \tilde\dbP_0 \circ (\tilde B)^{-1}$ is in $\overline\cP_W$. Denote
\beaa
E := \left\{\limsup_{t\downarrow 0} {\tilde B_{h}-\tilde B_0\over \sqrt{2h\ln\ln h^{-1}}} = 1\right\}.
\eeaa
Then $E\in \cF^{\tilde B}_{0+}$, and
$\tilde\dbP_0(E) = \dbP_0'(1) = {1\over 2}$.
\ep
}\end{eg}

\begin{eg} {\rm(Martingale measure without MRP)
\label{eg-MRP}
Let $\tilde\O:=(C[0,1])^2$, $(\tilde W,\tilde W')$
the canonical process, and $\tilde\dbP_0$ the
Wiener measure so that $\tilde W$ and $\tilde W'$
are independent Brownian motions under $\tilde\dbP_0$.
Let $\f: \dbR\to [0,1]$ be a measurable function, and
\beaa
\tilde B_t:=\int_0^t \tilde{\alpha}_s d\tilde W_s
&\mbox{where}&
\tilde{\alpha}_t:=[1+\f(\tilde W'_t)]^{1\over 2},~t\ge 0,
\eeaa
This induces the following probability measure
$\dbP$ on $\O$ with $d=1$,
\beaa
\dbP &:=& \tilde\dbP_0\circ \tilde B^{-1}.
\eeaa
Then $\dbP$ is a square integrable martingale measure
with $d\la B\ra_t\slash dt  \in [1,2]$, $\dbP$-almost surely.

We claim that $B$ has no MRP under $\dbP$. Indeed, if
$B$ has MRP under $\dbP$, then so does
$\tilde B$ under $\tilde\dbP_0$.
Let $\tilde\xi := \dbE^{\tilde\dbP_0}[\tilde W'_1|\cF^{\tilde B}_1]$.
Since $\tilde\xi \in \cF^{\tilde B}_1$ and is obviously
$\tilde\dbP_0$-square integrable, then there exists
$\tilde H^a\in \cH^2(\tilde\dbP_0, \dbF^{\tilde B})$ such that
\beaa
\tilde\xi = \dbE^{\tilde\dbP_0}[\tilde\xi] +\int_0^1
\tilde H^a_t d\tilde B_t =  \dbE^{\tilde\dbP_0}[\tilde\xi]
+\int_0^1 \tilde H^a_t \tilde\a_t d\tilde W_t, && \tilde\dbP_0-a.s..
\eeaa
Since $\tilde W$ and $\tilde W'$ are independent under
$\tilde\dbP_0$, we get $0 =  \dbE^{\tilde\dbP_0}[\tilde\xi \tilde W'_1]
= \dbE^{\tilde\dbP_0}[|\tilde\xi|^2]$. Then $\tilde\xi=0$,
$d\tilde\dbP_0$-almost surely, and thus
\bea\label{contradiction-noMRP}
\dbE^{\tilde\dbP_0}[\tilde W'_1 |\tilde B_1|^2]
&=&
\dbE^{\tilde\dbP_0}[\tilde\xi |\tilde B_1|^2]
\;=\;
0.
\eea
However, it follows from It\^o's formula,
together with the independence of $W$ and $W'$, that
\beaa
\dbE^{\tilde\dbP_0}[\tilde W'_1 |\tilde B_1|^2]
&=&
\dbE^{\tilde\dbP_0}\Big[\tilde W'_1 \int_0^1 2\tilde B_t \tilde\a_t d\tilde W_t\Big]
+
\dbE^{\tilde\dbP_0}\Big[\tilde W'_1\int_0^1 \tilde\a_t^2 dt\Big]\\
&=&
\dbE^{\tilde\dbP_0}\Big[\int_0^1 \tilde W'_t\big(1+\f(\tilde W'_t)\big) dt\Big]
=
\dbE^{\tilde\dbP_0}\Big\{\int_0^1 \tilde W'_t\f(\tilde W'_t)dt\Big\},
\eeaa
and we obtain a contradiction to \reff{contradiction-noMRP}
by observing that the latter expectation is non-zero for $\f(x) := \1_{\dbR_+}(x)$.
\ep
}\end{eg}

We note that, however, we are not able to find a
good example such that  $a\in {\cA}_W$ (so that
\reff{SDE0} has unique weak solution) but $B$ has
no MRP under $\dbP^a$ (and consequently \reff{SDE0} has no strong solution).

\subsection{Some technical proofs}

\no{\it Proof of Lemma \ref{lem-version}.}
The uniqueness is obvious. We now prove the existence.

\no (i) Assume $X$ is {\cad}, $\dbP$-almost surely.
Let $E_0:=\{\o: X_\cd(\o)$ is not \cad $\}$.
For each $r\in \dbQ\cap(0,\infty)$, there exists
$\tilde X_r\in \cF^+_{r}$ such that
$E_r:=\{\tilde X_r \neq X_r\}\in \cN^\dbP(\cF_\infty)$.
Let $E:= E_0\cup(\cup_r E_r)$. Then $\dbP(E)=0$.
For integers $n\ge 1$ $k\ge 0$, set $t^n_k:=k/n$, and define
\beaa
X^n_t := \tilde X_{t^n_{k+1}}~\mbox{for}~t\in\big(t^n_{k},t^n_{k+1}\big],
&\mbox{and}&
\tilde X:= (\limsup_{n\to\infty} X^n )\1_{\{\limsup_{n\to\infty} X^n\in\dbR\}}.
\eeaa
Then for any $t\in (t^n_k, t^n_{k+1}]$,
$X^n_t\in \cF^+_{t^n_{k+1}}$ and
$X^n|_{[0,t]}\in \cB([0,t])\times \cF^+_{t^n_{k+1}}$.
Since $\dbF^+$ is right continuous, we get
$\tilde X_t\in \cF^+_{t}$ and $\tilde X|_{[0,t]}\in \cB([0,t])\times \cF^+_{t}$.
That is, $\tilde X\in\dbF^+$.
Moreover, for any $\o\notin E$ and $n\ge 1$, if
$t\in (t^n_k, t^n_{k+1}]$, we get
$$
\lim_{n\to\infty} X^n_t(\o) =
\lim_{n\to\infty}\tilde X_{t^n_{k+1}}(\o) =  \lim_{n\to\infty}X_{t^n_{k+1}}(\o) = X_t(\o).
$$
So $\{\o: ~\mbox{there exists}~t\ge 0 ~\mbox{such that}
~~ \tilde X_t(\o) \neq X_t(\o)\}\subset E$.
Then, $\tilde X$ is $\dbP$-indistinguishable
from $X$ and thus $\tilde X$ also has {\cad} paths, $\dbP$-almost surely.

\ms
\no (ii) Assume $X$ is $\overline\dbF^\dbP$-progressively
measurable and is bounded. Let $Y_t:= \int_0^t X_s ds$.
Then $Y$ is continuous.
By (i), there exists $\dbF^+$-progressively measurable
continuous process $\tilde Y$ such that $\tilde Y$ and
$Y$ are $\dbP$-indistinguishable. Let
$E_0:=\{\mbox{there exists}~ t\ge 0 ~ \mbox{such that}
~\tilde Y_t \neq Y_t\}$, then $\dbP(E_0)=0$ and
$\tilde Y_\cd(\o)$ is continuous for each $\o\notin E_0$. Define,
\beaa
X^n_t := n[\tilde Y_{t}-\tilde Y_{t-\frac{1}{n}}];\q \tilde
X:= (\limsup_{n\to\infty} X^n )\1_{\{\limsup_{n\to\infty} X^n\in\dbR\}}
&\mbox{for}&
n\ge 1.
\eeaa
As in (i), we see $\tilde X\in\dbF^+$. Moreover, for each
$\o\notin E_0$, $X^n_t(\o) = n\int_{t-\frac{1}{n}}^t X_s(\o)ds$.
Then $\tilde X_\cd(\o)=X_\cd(\o), dt$-almost surely.
Therefore, $\tilde X=X$, $\dbP$-almost surely.

\ms

\no (iii) For general $\overline\dbF^\dbP$-progressively
measurable $X$, let $X^m_t:= (-m)\vee (X\wedge m)$,
for any $m\ge 1$. By (ii), $X^m$ has an $\dbF^+$-adapted
modification $\tilde X^m$. Then obviously the following process $\tilde X$ satisfies all the requirements:
$\tilde X:= (\limsup_{m\to\infty} \tilde X^m )
\1_{\{\limsup_{m\to\infty} \tilde X^m\in\dbR\}}$.
\ep

\vspace{5mm}

To prove Example \ref{eg-Piecewise constant}, we need a simple lemma.

\begin{lem}
\label{lem-weaktau}
Let $\t$ be an $\dbF$-stopping time and
$X$ is an $\dbF$-progressively measurable process.
Then $\t(X_\cd)$ is also an $\dbF-$stopping time.

Moreover, if $Y$ is $\dbF$-progressively measurable
and $Y_t = X_t$ for all $t\le \t(X_\cd)$, then $\t(Y_\cd) = \t(X_\cd)$.
\end{lem}

\proof Since $\t$ is an $\dbF$-stopping time,
we have $\{\t(X_\cd) \le t\} \in \cF^X_t$ for all $t\ge 0$.
Moreover, since $X$ is $\dbF$-progressively measurable,
we know $\cF^X_t \subset \cF^B_t$.
Then $\{\t(X_\cd)\le t\}\in \cF^B_t$ and thus $\t(X_\cd)$ is
an $\dbF-$stopping time.

Now assume $Y_t = X_t$ for all $t\le \t(X_\cd)$.
For any $t\ge 0$, on $\{\t(X_\cd) = t\}$, we have
$Y_s = X_s$ for all $s\le t$. Since $\{\t(X_\cd) = t\}\in \cF^X_t$ and
by definition $\cF^X_t = \si(X_s, s\le t\}$,
then $\t(Y_\cd)=t$ on the event $\{\t(X_\cd) = t\}$. Therefore, $\t(Y_\cd)=\t(X_\cd)$.
\ep

\ms

\no{\it Proof of Example \ref{eg-Piecewise constant}.}
Without loss of generality we prove only that
\reff{SDE0} on $\dbR_+$ with $X_0=0$ has a unique strong solution.
In this case the stochastic differential equation becomes
\beaa
dX_t = \sum_{n=0}^\infty a_n(X_\cd) \1_{[\t_n(X_\cd),
\t_{n+1}(X_\cd))} dB_t, &t\ge 0,& \dbP_0-a.s..
\eeaa
We prove the result by induction on $n$. Let $X^0$ be the solution to SDE:
\beaa
X^0_t = \int_0^t a_0^{{1\slash 2}}(X^0_\cd) dB_s,
&t\ge 0&, \dbP_0-\mbox{almost surely}
\eeaa
Note that $a_0$ is a constant, thus
$X^0_t = a_0^{1\over 2}B_t$ and is unique.
Denote $\tilde \t_0 := 0$ and $\tilde\t_1 := \t_1(X^0_\cd)$.
By Lemma \ref{lem-weaktau}, $\tilde \t_1$ is an $\dbF-$stopping time.
Now let $X^1_t := X^0_t$ for $t\le \tilde\t_1$, and
\beaa
X^1_t = X^0_{\tilde\t_1} + \int_{\tilde\t_1}^t a_1^{{1\slash 2}}
(X^1_\cd) dB_s, &t\ge \tilde\t_1,& \dbP_0-\mbox{a.s.}
\eeaa
Note that $a_1\in\cF_{\t_1}$, that is, for any $y\in\dbR$ and
$t\ge 0$, $\{a_1(B_\cd)\le y, \t_1(B_\cd)\le t\}\in\cF_t$.
Thus, for any $\bx, \tilde \bx\in C(\dbR_+, \dbR^d)$, if
$\bx_s = \tilde \bx_s, 0\le s\le t$, then $a_1(\bx)\1_{\{\t_1(\bx)\le t\}}
= a_1(\tilde \bx)\1_{\{\t_1(\tilde \bx)\le t\}}$. In particular,
noting that $\t_1(X^1_\cd) = \t_1(X^0_\cd) = \tilde\t_1$,
for each $\o$ by choosing $t=\tilde\t$
we obtain that $a_1(X^1_\cd) = a_1(X^0_\cd)$.
Thus $X^1_t = X^0_{\tilde \t_1} + a_1(X^0_\cd)[B_t-B_{\tilde\t_1}]$,
$t\ge\tilde\t_1$, and is unique.
Now repeat the procedure for $n=1,2,\cds$
we obtain the unique strong solution
$X$ in $[0, \tilde\t_\infty)$, where $\tilde\t_{\infty} := \lim_{n\to\infty}\t_n(X_\cd)$.
Since $a$ is bounded, it is obvious that $X_{\tilde\t_\infty}:=\lim_{t\uparrow \tilde\t_\infty}
X_t$ exists $\dbP_0$-almost surely.
Then, by setting $X_t := X_{\tilde\t_\infty}$ for $t\in (\tilde\t_{\infty},\infty)$
we complete the construction.
\ep

\bs

\no{\it Proof of Lemma \ref{lem-Astructure}.}
Let $a$ be given as in \reff{Aa} and $\t\in\cT$ be fixed.
First, since $\{E^n_i, i\ge 1\}$ is a partition of $\O$, then for any $n\ge 0$,
\beaa
\left\{\cap_{j=0}^n E^j_{i_j}, ~~(i_j)_{0\le j\le n}\in\dbN^{n+1}\right\}
~~\mbox{also form a partition of}~\O.
\eeaa
Next, assume $\t_n$ takes values $t^n_k$
(possibly including the value $\infty$), $k\ge 1$.
Then $\{\{\t_n = t^n_k\}, k\ge 1\}$ form a partition of $\O$.
Similarly we have, for any $n\ge 0$,
\beaa
\left\{\cap_{j=0}^{n+1}\{\t_j = t^j_{k_j}\},~(k_j)_{0\le j\le n+1}
\in\dbN^{n+2}\right\}
&\mbox{form a partition of}&
\O.
\eeaa
These in turn form another partition of $\O$ given by,
\bea
\label{partition}
\left\{\Big[\cap_{j=0}^n \big(E^j_{i_j}\cap \{\t_j = t^j_{k_j}\}\big)\Big]\bigcap
\{\t_{n+1}=t^{n+1}_{k_{n+1}}\}, ~~(i_j, k_j)_{0\le j\le n}
\in\dbN^{2(n+1)},~ k_{n+1}\in \dbN\right\}.
\eea
Denote by $\cI$ the family of all finite sequence of
indexes $I:=(i_j, k_j)_{0\le j\le n}$ for some $n$
such that $0=t^0_{k_0} < \cds< t^{n}_{k_{n}}<\infty$.
Then $\cI$ is countable. For each $I\in\cI$, denote by
$|I|$ the corresponding $n$, and define
\beaa
E_I &:=& \left(\cap_{j=0}^{|I|} \left[E^j_{i_j}\cap \{\t_j = t^j_{k_j}
\le \t\}\right]\right)\bigcap\left(\{\t_{|I|+1}>\t\}\cup \{\t_{|I|+1}=\t=\infty\}\right),\\
\tilde\t &:=& \sum_{I\in\cI} \t_{|I|+1}\1_{E_I},
~~\mbox{and}~~
a_I \;:=\; \sum_{j=0}^{|I|-1} a^j_{i_j}\1_{[t^j_{k_j}, t^{j+1}_{k_{j+1}})}
+ a^{|I|}_{i_{|I|}}\1_{[t^{|I|}_{k_{|I|}}, \infty)}.
\eeaa
It is clear that $E_I$ is $\cF_{\t}-$measurable.
Then, in view of the concatenation property of $\cA_0$,
$a_I\in\cA_0$.
In light of \reff{partition},  we see that $\{E_I, I\in\cI\}$
are disjoint. Moreover, since $\t_n = \infty$ for $n$ large enough,
we know $\{E_I, I\in\cI\}$ form a partition of $\O$.
Then $\tilde\t$ is an $\dbF-$stopping time and either
$\tilde\t>\t$ or $\tilde\t=\t=\infty$.
We now show that
\bea
\label{aI}
a_t = \sum_{I\in\cI} a_I(t)\1_{E_I} &\mbox{for all}& t<\tilde\t.
\eea
In fact, for each $I =(i_j, k_j)_{0\le j\le n}\in\cI$, $\o\in E_I$, and $t<\tilde\t(\o)$, we have $\t_j(\o) = t^j_{k_j}\le \t(\o)$ for $j\le n$ and $\t_{n+1}(\o) = \tilde\t(\o)>t$. Let $j_0 = j_0(t,\o)\le n$ be such that $\t_{j_0}(\o) \le t < \t_{j_0+1}(\o)$. Then $\1_{[\t_{j_0}(\o), \t_{j_0+1}(\o))}(t) = 1$ and $\1_{[\t_j(\o), \t_{j+1}(\o))}(t) = 0$ for $j\neq j_0$, and thus
\beaa
a_t(\o) &=& \sum_{j=0}^\infty \sum_{i=1}^\infty a^j_i(t,\o) \1_{E^j_i}(\o)\1_{[\t_j(\o), \t_{j+1}(\o))}(t)=\sum_{i=1}^\infty a^{j_0}_i(t,\o) \1_{E^{j_0}_i}(\o) = a^{j_0}_{i_{j_0}}(t,\o),
\eeaa
where the last equality is due to the fact that $\o\in E_I \subset E^{j_0}_{i_{j_0}}$ and that $\{E^{j_0}_i, i\ge 1\}$ is a partition of $\O$. On the other hand, by the definition of $a_I$, it is also straightforward to check that $a_I(t,\o) = a^{j_0}_{i_{j_0}}(t,\o)$. This proves \reff{aI}. Now since $\cI$ is countable, by numerating the elements of $\cI$ we prove the lemma.

Finally, we should point out that, if $\t=\t_n$, then we can choose $\tilde\t = \t_{n+1}$.
\ep

\end{document}